\documentclass{article}

\usepackage{amsmath}
\usepackage{amsthm}
\usepackage{amssymb}

\title{Deformation of locally free actions\\
 and the leafwise cohomology}
\author{Masayuki Asaoka\\
Kyoto University, Japan}

\def\ra{{\rightarrow}}
\def\Id{{\text{\rm Id}}}
\DeclareMathOperator{\Diff}{Diff}
\DeclareMathOperator{\Img}{Im}
\DeclareMathOperator{\Ker}{Ker}
\DeclareMathOperator{\Hom}{Hom}

\def\RR{\mathbb R}
\def\TT{\mathbb T}
\def\CC{\mathbb C}
\def\ZZ{\mathbb Z}
\def\HH{\mathbb H}

\def\cA{{\cal A}}
\def\cF{{\cal F}}
\def\cO{{\cal O}}
\def\cU{{\cal U}}
\def\cT{{\cal T}}
\def\st{{\;|\;}}

\def\frX{{\mathfrak{X}}}
\def\frg{{\mathfrak{g}}}
\def\frh{{\mathfrak{h}}}
\def\frga{{\mathfrak{ga}}}

\def\del{\partial}
\def\--{{\setminus}}
\def\bs{{\backslash}}

\def\hsp{{\hspace{3mm}}}

\theoremstyle{plain}
\newtheorem{thm}{Theorem}[section]
\newtheorem{prop}[thm]{Proposition}
\newtheorem{lemma}[thm]{Lemma}
\newtheorem{cor}[thm]{Corollary}

\theoremstyle{definition}

\newtheorem{conj}[thm]{Conjecture}
\newtheorem{expl}[thm]{Example}
\newtheorem{excs}[thm]{Exercise}

\theoremstyle{remark}

\begin{document}
\maketitle

\section{Introduction}

This is a note of the author's lectures
 at the ``Advanced Course on Foliations''
 in the research program ``Foliations'',
 which was held at the Centre de Recerca Mathem\`atica
 in May of 2010.
In this note, we discuss the relationship
 between deformations of actions of Lie groups
 and the leafwise cohomology of the orbit foliation.

In early 1960's,
Palais \cite{Pa61} proved the local rigidity of
smooth actions of compact groups.
Hence, such actions have no non-trivial deformations.
In contrast to compact groups,
Any known $\RR$-actions ({\it i.e.,} flows)
 are not locally rigid except trivially rigid ones.
Moreover, many of $\RR$-actions change
 the topological structure of their orbits under perturbation.
Their bifurcation is an important issue in the theory of dynamical systems.

In the last two decades, it has been found that
 there exist locally rigid actions of higher dimensional Lie groups,
 and the rigidity theory of locally free actions
 have been rapidly developed.
The reader can find examples of locally rigid or
parameter rigid actions
in many papers \cite{As09,Da07,DK05,DK07,DKpre,dS93,dS07,EN93,
Gh85,KS94,KS97,KY97,MM03,Mi07,Wa10,Wa10-2},
some of which we will discuss in this article.

Rigidity problem can be regarded as a special case of deformation problem.
In many situations,
 the deformation space of a geometric structure is described 
 by a system of non-linear partial differential equations.
Its linearization defines a cochain complex,
 so called {\it the deformation complex},
 and the space of infinitesimal deformations is
 identified with the first cohomology of this complex.
For locally free actions of Lie groups,
 the deformation complex is realized as
 the (twisted) leafwise de Rham complex of the orbit foliation.

The reader may wish to develop a general deformation theory
 of locally free actions in terms of the deformation complex,
 like the deformation theory of complex manifolds
 founded by Kodaira and Spencer.
However, the leafwise de Rham complex is not elliptic,
 and this causes two difficulties to develop a fine theory.
First, the leafwise cohomology groups are infinite dimensional in general,
 and they are hard to compute.
Second, we need to apply the implicit function theorem for maps
 between Fr\'echet spaces rather than Banach spaces
 because of ``loss of derivative''.
It requires tameness of splitting of the deformation complex,
 which is hard to show.
So, we will focus on the techniques to overcome these difficulties
 for several explicit examples, instead of developing a general theory.

The main tools for computation of the leafwise cohomology 
 are Fourier analysis, representation theory,
 and Mayer-Vietoris argument developed by El Kacimi Aloui
 and Tihami.
Matsumoto and Mitsumatsu also developed a technique based on
 ergodic theory of hyperbolic dynamics.
We will discuss these techniques
 in Section \ref{sec:leafwise cohomology}.

For several actions,
 deformation problem can be reduced to a linear one
 without help of any implicit function theorem,
 and hence, we can avoid tame estimate of the splitting.
In Section \ref{sec:orbit preserving},
 we will see how to reduce rigidity problem of such actions
 to (almost) vanishing of the first cohomology
 of the leafwise cohomology.
The first case is parameter deformation of abelian actions.
We will see the problem is a linear one in this case.
In fact, the deformation space can be naturally identified
 with the space of infinitesimal deformations.
The second case is parameter rigidity of solvable actions.
Although the problem itself is not linear in this case,
 we can decompose it into the solvability of linear equations
 for several examples.

For general cases,
 the deformation problem cannot be reduced to a linear one directly.
 One way to describe the deformation space
 is to apply Hamilton's implicit function theorem.
As mentioned above,
 it needs a tame estimate on solutions of partial differential equations
 and it is difficult to establish it in general.
But, there are a few examples to which we can apply the theorem.
Another way is to use the theory of hyperbolic dynamics.
We will give a brief discussion on
 these techniques in Section \ref{sec:general}.

The author recommend the readers
 to read the survey papers \cite{BE} and \cite{Ma09}.
The former contains a nice exposition
 on applications of Hamilton's implicit function theorem
 to rigidity problems of foliations.
The second is a survey on the parameter rigidity problem,
 which is one of the sources of the author's lectures at 
 the Centre de Recerca Matem\`atica.
\bigskip

To end the introduction,
 the author would like to thank the organizers
 of the research program ``Foliations'' at the CRM
 for inviting me to give these lectures in the program,
 and the staff of the CRM for their warm hospitality.
The author is also grateful to Marcel Nicolau
 and Nathan dos Santos for many suggestions to improve this note.

\section{Locally free actions and their deformations}
\label{sec:locally free actions}

In this section, we define locally free actions
and their infinitesimal correspondent.
We also introduce the notion of deformation of actions
and several concepts of finiteness of
codimension of the conjugacy classes of an action
in the space of locally free actions.

\subsection{Locally free actions}
In this note,
 we will work in the $C^\infty$-category.
So, the term ``smooth'' means ``$C^\infty$''
 and all diffeomorphisms are of class $C^\infty$.
All manifolds and Lie groups in this note will be connected.
For manifolds $M_1$ and $M_2$,
 we denote the space of smooth maps from $M_1$ to $M_2$
 by $C^\infty(M_1,M_2)$.
It is endowed with the $C^\infty$ compact-open topology.
By $\cF(x)$, we denote the leaf of a foliation $\cF$
 which contains a point $x$.

Let $G$ be a Lie group and $M$ a manifold.
We denote the unit element of $G$ by $1_G$
 and the identity map of $M$ by $\Id_M$.
We say a smooth map $\rho:M \times G \ra M$ is
 {\it a (smooth right) action} if
\begin{enumerate}
 \item $\rho(x,1_G)=x$ for any $x \in M$, and
 \item $\rho(x,gh)=\rho(\rho(x,g),h)$
 for any $x \in M$ and $g,h \in G$.
\end{enumerate}
For $\rho \in C^\infty(M \times G, M)$ and $g \in G$,
 we define a map $\rho^g:M \ra M$ by $\rho^g(x)=\rho(x,g)$.
Then, $\rho$ is an action if and only if
 the map $g \mapsto \rho^g$ is an anti-homomorphism from $G$ into
 the group $\Diff^\infty(M)$ of diffeomorphisms of $M$.
By $\cA(M,G)$, we denote the subset of
 $C^\infty(M \times G,M)$ that consists of actions of $G$.
It is a closed subspace of $C^\infty(M \times G,M)$.
For $\rho \in \cA(M,G)$ and $x \in M$, the set
\begin{equation*}
 \cO_\rho(x)=\{\rho^g(x) \st g \in G\}
\end{equation*}
 is called {\it the $\rho$-orbit} of $x$.

\begin{expl}
$\cA(M,G)$ is non-empty for any $M$ and $G$.
In fact, it contains {\it the trivial action} $\rho_{triv}$,
 which is defined by $\rho_{triv}(x,g)=x$.
For any $x \in M$, $\cO_{\rho_{triv}}(x)=\{x\}$.
\end{expl}

Let us introduce an infinitesimal description of actions.
By $\frX(M)$, we denote the Lie algebra of smooth vector fields on $M$.
Let $\frg$ be the Lie algebra of $G$
 and $\Hom(\frg,\frX(M))$ be the space of Lie algebra homomorphisms
 from $\frg$ to $\frX(M)$.
In this note, we identify $\frg$
 with the subspace of $\frX(G)$ consisting of
 vector fields invariant under left translations.
Each action $\rho \in \cA(M,G)$ associates
 {\it the infinitesimal action } $I_\rho:\frg \ra \frX(M)$ by
\begin{equation*}
\left.I_\rho(\xi)(x)=\frac{d}{dt}\rho(x,\exp t \xi)\right|_{t=0}.
\end{equation*}
\begin{prop}
\label{prop:I is homomorphism} 
$I_\rho$ is a Lie algebra homomorphism from $\frg$ to $\frX(M)$.
\end{prop}
\begin{proof}
By ${\cal L}_X Y$, we denote the Lie derivative of a vector field $Y$ 
 by another vector field $X$.
Take $\xi,\eta \in \frg$ and $x \in M$.
Then,
\begin{align*}
[I_\rho(\xi),I_\rho(\eta)](x)
 & =({\cal L}_{I_\rho(\xi)}I_\rho(\eta))(x)\\
 & =\lim_{t \ra 0}\frac{1}{t}
 \left\{D\rho^{\exp(-t\xi)}(I_\rho(\eta)(\rho^{\exp(t\xi)}(x)))
 -I_\rho(\eta)(x) \right\}\\
 & = \left.\frac{d}{dt}\left\{
 \left. \frac{d}{ds}
 (\rho^{\exp(-t \xi)} \circ \rho^{\exp(s \eta)} \circ \rho^{\exp(t \xi)})(x)
 \right|_{s=0}
 \right\} \right|_{t=0}\\
 & = \left.\frac{d}{dt}\left\{
 \left. \frac{d}{ds}
 \rho(x,\exp(t \xi)\exp(s \eta)\exp(-t \xi))
 \right|_{s=0}
 \right\} \right|_{t=0}\\
 & = \left.\frac{d}{dt}
 \rho(x,\exp(t \text{Ad}_{\exp(t\xi)}\eta))
  \right|_{t=0}\\
 & = I_\rho([\xi,\eta])(x).
\end{align*}
\end{proof}

\begin{prop}
\label{prop:infinitesimal action} 
Two actions $\rho_1, \rho_2 \in \cA(M,G)$
 coincide if $I_{\rho_1}=I_{\rho_2}$.
If $G$ is simply connected and $M$ is closed,
 then  any $I \in \Hom(\frg,\frX(M))$ is
 the infinitesimal action associated with some action in $\cA(M,G)$.
\end{prop}
\begin{proof}
Take $\rho_1,\rho_2 \in \cA(M,G)$.
The curve $t \mapsto \rho_i(x,\exp(t\xi))$ is an integral curve
 of the vector field $I_{\rho_i}(\xi)$
 for any $i=1,2$, $x \in M$ and $\xi \in \frg$.
If $I_{\rho_1}=I_{\rho_2}$, then
 the uniqueness of the integral curve implies
 that $\rho_1(x,\exp(t\xi))=\rho_2(x,\exp(t\xi))$
 for any $x \in M$, $t \in \RR$, and $\xi \in \frg$.
Since the union of one-parameter subgroups of $G$ generates $G$,
 we have $\rho_1=\rho_2$.
 
Suppose that $G$ is simply connected and $M$ is a closed manifold.
Let $E$ be a subbundle of $T(M \times G)$ given by
\begin{equation*}
E(x,g)=\{(I(\xi)(x),\xi(g)) \in T_{(x,g)} (M \times G) \st \xi \in \frg\}.
\end{equation*}
For any $\xi, \xi' \in \frg$, we have
\begin{align*}
[(I(\xi),\xi), (I(\xi'),\xi')] 
 =([I(\xi),I(\xi')],[\xi,\xi']) = (I([\xi,\xi']),[\xi,\xi']).
\end{align*}
By Frobenius' theorem, the subbundle $E$ is integrable.
Let $\cF$ be the foliation on $M \times G$ generated by $E$.
The space $M \times G$ admits a left action of $G$
 by $g \cdot (x,g')=(x,gg')$.
The subbundle $E$ is invariant under this action.
Hence, we have $g \cdot \cF(x,g')=\cF(x,gg')$.
Since $M$ is compact, $G$ is simply connected,
 and the foliation $\cF$ is transverse to the
 natural fibration $\pi:M \times G \ra G$,
 the restriction of $\pi$ to each leaf of $\cF$
 is a diffeomorphism onto $G$.
So, we can define a smooth map $\rho:M \times G \ra M$ 
 so that $\cF(x,1_G) \cap \pi^{-1}(g)=\{(\rho^g(x),g)\}$.
Take $x \in M$ and $g,h \in G$.
Then, $(\rho^g \circ \rho^{h}(x),g)$ is contained in $\cF(\rho^{h}(x),1_G)$.
Applying $h$ from left,  we see that
 $(\rho^g \circ \rho^{h}(x),hg)$ is an element of $\cF(\rho^{h}(x),h)$.
Since $\cF(\rho^{h}(x),h)=\cF(x,1_G)$
 and $\{(\rho^{hg}(x), hg)\}=\cF(x,1_G) \cap \pi^{-1}(hg)$
 by the definition of $\rho$,
 we have $\rho^g \circ \rho^{h}(x)=\rho^{hg}(x)$.
Therefore, $\rho$ is a right action of $G$.
Now, it is easy to check that $I_\rho = I$.
\end{proof}

We say that an action $\rho \in \cA(M,G)$ is {\it locally free}
 if the isotropy group $\{g \in G \st \rho^g(x)=x\}$
 is a discrete subgroup of $G$ for any $x \in M$.
By $\cA_{LF}(M,G)$,
 we denote the set of locally free actions of $G$ on $M$.
Of course, the trivial action is not locally free
 unless $M$ is zero-dimensional.
The following is a list of basic examples of locally free actions.
\begin{expl}
[Flows]
A locally free $\RR$-action is just a smooth flow
 with no stationary points.
Remark that $\cA_{LF}(M,\RR)$ is empty
 if $M$ is a closed manifold with non-zero Euler characteristic.
\end{expl}
\begin{expl}
[The standard action]
Let $G$ be a Lie group,
 and $\Gamma$, $H$ be closed subgroups of $G$.
{\it The standard $H$-action on $\Gamma \bs G$}
 is the action $\rho_\Gamma \in \cA(\Gamma \bs G,H)$ defined by
 $\rho(\Gamma g,h)=\Gamma(gh)$.
The action $\rho$ is locally free if and only if
 $g^{-1}\Gamma g \cap H$ is a discrete subgroup of $H$
 for any $g \in G$.
In particular,
 if $\Gamma$ itself is a discrete subgroup of $G$,
 then $\rho$ is locally free.
\end{expl}
\begin{expl}
[The suspension construction] 
Let $M$ be a manifold and $G$ be a Lie group.
Take a discrete subgroup $\Gamma$ of $G$,
 a closed subgroup $H$ of $G$,
 and a left action $\sigma:\Gamma \times M  \ra M$.
We put $M \times_\sigma G
 = M \times G / (x,g) \sim (\sigma(\gamma, x),\gamma g)$.
Then, $M \times_\sigma G$ is an $M$-bundle over $\Gamma \bs G$.
We define a locally free action $\rho$ of $H$ on $M \times_\sigma G$
 by $\rho([x,g],h)=[x,gh]$.
\end{expl}

We say a homomorphism $I:\frg \ra \frX(M)$
 is {\it non-singular} if $I(\xi)(x) \neq 0$
 for any $\xi \in \frg \-- \{0\}$ and $x \in M$.
\begin{prop}
\label{prop:locally free} 
An action $\rho \in \cA(M,G)$ is locally free 
 if and only if $I_\rho$ is non-singular.
\end{prop}
\begin{cor}
\label{cor:locally free} 
For any $\rho \in \cA_{LF}(M,G)$,
 the orbits of $\rho$ form a smooth foliation.
If the manifold $M$ is closed,
 then the map $\rho(x,\cdot):G \ra \cO(x,\rho)$ is a covering
 for any $x \in M$,
 where $\cO(x,\rho)$ is endowed with the leaf topology.
\end{cor}
Proofs of the proposition and the corollary are easy
 and left to the reader.
If $M$ is closed, the set of non-singular homomorphisms
 is an open subset of $\Hom(\frg,\frX(M))$.
Hence, $\cA_{LF}(M,G)$ is an open subset of $\cA(M,G)$ in this case.

Let $\cF$ be a foliation on a manifold $M$.
We denote the tangent bundle of $\cF$ by $T\cF$,
 and the subalgebra of $\frX(M)$ consisting of
 vector fields tangent to $\cF$ by $\frX(\cF)$.
Let $\cA_{LF}(\cF,G)$ be the set of locally free actions 
 of a Lie group $G$ whose orbit foliation is $\cF$.
The subspace $\cA_{LF}(\cF,G)$ of $\cA_{LF}(M,G)$ is closed
 and it consists of actions $\rho$
 such that $I_\rho$ is an element of $\Hom(\frg,\frX(\cF))$.

\subsection{Rigidity and deformations of actions}
We say that two actions $\rho_1 \in \cA(M_1,G)$
 and $\rho_2 \in \cA(M_2,G)$ on manifolds $M_1$ and $M_2$
 are {\it ($C^\infty$-)conjugate}
 (and write $\rho_1 \simeq \rho_2$) if
 there exists a diffeomorphism $h:M_1 \ra M_2$
 and an automorphism $\Theta$ of $G$ such that
 $\rho_2^{\Theta(g)} \circ h= h \circ \rho_1^g$
 for any $g \in G$.
For a given foliation $\cF$ on $M$,
 let $\Diff(\cF)$ be the set of diffeomorphisms of $M$
 which preserve each leaf of $\cF$,
 and $\Diff_0(\cF)$ be its arc-wise connected component
 that contains $\Id_M$.
We say that two actions $\rho_1,\rho_2 \in \cA_{LF}(\cF,G)$
 are {\it ($C^\infty$) parameter-equivalent}
 (and write $\rho_1 \equiv \rho_2$) if
 they are conjugate by a pair $(h,\Theta)$
 such that $h$ is an element of $\Diff_0(\cF)$.
It is easy to see that
 conjugacy and the parameter-equivalence are equivalence relations.

The ultimate goal
 is the classification of actions in $\cA_{LF}(M,G)$ or $\cA_{LF}(\cF,G)$
 up to conjugacy, or parameter-equivalence
 for given $G$ and $M$, or $\cF$.
The simplest case is that $\cA_{LF}(M,G)$ or $\cA_{LF}(\cF,G)$ consists of
 only one equivalence class.
We say that an action $\rho_0$ in $\cA_{LF}(M,G)$
 is ($C^\infty$-){\it rigid}
 if any action in $\cA_{LF}(M,G)$ is conjugate to
 $\rho_0$.
We say that an action $\rho_0$ whose orbit foliation is $\cF$
 is ($C^\infty$-){\it parameter rigid}
 if any action in $\cA_{LF}(\cF,G)$ is parameter-equivalent to $\rho_0$.

It is useful to introduce a local version of rigidity.
We say $\rho_0$ is {\it locally rigid}
 if there exists a neighborhood $\cU$ of $\rho_0$
 such that any action in $\cU$ is conjugate
 to $\rho_0$.
We also say $\rho_0$ is {\it locally parameter rigid}
 if there exists a neighborhood $\cU$ of $\rho_0$
 in $\cA(\cF,G)$
 such that any action in $\cU$ is parameter-equivalent
 to $\rho_0$\footnote{
There exists an action which is
 locally parameter rigid but not parameter rigid.
For $k \in \ZZ$, let $\rho_k$ be a right action of $S^1=\RR/\ZZ$
 on $S^1$ by $\rho_k^t(s)=s+kt$.
It is easy to see that $\rho_1$ is locally parameter rigid.
Of course,
 all the orbits of $\rho_k$ coincides with $S^1$ for any $k \geq 1$.
However, $\rho_k$ is  parameter equivalent to $\rho_1$
 if and only if $|k|=1$
 since the mapping degree of $\rho_k(s,\cdot)$ is $k$.
So, $\rho_1$ is locally parameter rigid, but is not parameter rigid.

It is unknown whether any locally parameter rigid locally free action
 of contractible Lie group on a closed manifold
 is parameter rigid or not.}.
As we mentioned in the introduction,
 local rigidity for compact group actions was settled in early 1960's.
\begin{thm}
[Palais \cite{Pa61}]
Any action of a compact group on a closed manifold 
 is locally rigid.
\end{thm}

As we see later,
 many actions of non-compact groups are not locally rigid.
So, it is natural to introduce the concept of deformation of actions.
We say that a family $(\rho_\mu)_{\mu \in \Delta}$
 of elements of $\cA(M,G)$ parametrized by a manifold $\Delta$
 is {\it a $C^\infty$ family} 
 if the map $\bar{\rho}:(x,g,\mu) \mapsto \rho_\mu(x,g)$
 is a smooth map.
By $\cA_{LF}(M,G;\Delta)$,
 we denote the set of $C^\infty$ family 
 of actions in $\cA_{LF}(M,G)$ parametrized by $\Delta$.
Under the identification with
 $(\rho_\mu)_{\mu \in \Delta}$ and $\bar{\rho}$,
 the topology of $C^\infty(M \times G \times \Delta, M)$
 induces a topology on $\cA_{LF}(M,G;\Delta)$.
We say that $(\rho_\mu)_{\mu \in \Delta}$ is
 {\it a (finite dimensional) deformation} of $\rho \in \cA(M,G)$
 if $\Delta$ is an open neighborhood of $0$
 in a finite dimensional vector space
 and $\rho_0=\rho$.

In several cases, actions are not locally rigid,
 but their conjugacy class is ``of finite codimension'' in $\cA_{LF}(M,G)$.
Here, we formulate two types of finiteness of codimension.
Let $(\rho_\mu)_{\mu \in \Delta} \in \cA_{LF}(M,G;\Delta)$
 be a deformation of $\rho$.
We say that $(\rho_\mu)_{\mu \in \Delta}$ is {\it locally complete} 
 if there exists a neighborhood $\cU$ of $\rho$ in $\cA_{LF}(M,G)$ 
 such that any action in $\cU$
 is conjugate to $\rho_\mu$ for some $\mu \in \Delta$.
We also say that $(\rho_\mu)_{\mu \in \Delta}$
 is {\it locally transverse}\footnote{
 This terminology is not common. 
Any suggestion of a better terminology is welcome.}
 if any $C^\infty$ family in $\cA_{LF}(M,G;\Delta)$
 sufficiently close to $(\rho_\mu)_{\mu \in \Delta}$
 contains an action conjugate to $\rho$.
Roughly speaking,
 the local completeness means
 that the quotient space $\cA(M,G)/\simeq$ is locally finite dimensional
 at the conjugacy class of $\rho$.
The local transversality means
 the family $(\rho_\mu)_{\mu \in \Delta}$ is transverse to
 the conjugacy class of $\rho$ at $\mu=0$.

We define analogous concepts for actions in $\cA_{LF}(\cF,G)$.
Let $\cF$ be a foliation on a manifold $M$.
We say that $(\rho_\mu)_{\mu \in \Delta} \in \cA_{LF}(M,G;\Delta)$
 {\it preserves} $\cF$
 if all $\rho_\mu$'s are actions in $\cA_{LF}(\cF,G)$.
By $\cA_{LF}(\cF,G;\Delta)$,
 we denote the subset of $\cA_{LF}(M,G;\Delta)$
 that consists of families preserving $\cF$.
We call a deformation in $\cA_{LF}(\cF,G;\Delta)$
 {\it a parameter deformation}.
Let $(\rho_\mu)_{\mu \in \Delta} \in \cA_{LF}(\cF,G;\Delta)$
 be a parameter deformation of an action $\rho$.
We say that $(\rho_\mu)_{\mu \in \Delta}$
 is {\it locally complete in $\cA_{LF}(\cF,G)$}
 if there exists a neighborhood $\cU$ of $\rho$ in $\cA_{LF}(\cF,G)$ 
 such that any action in $\cU$
 is parameter-equivalent to $\rho_\mu$ for some $\mu \in \Delta$.
We also say that
 $(\rho_\mu)_{\mu \in \Delta} \in \cA_{LF}(\cF,G;\Delta)$
 is {\it locally transverse in $\cA_{LF}(\cF,G)$}
 if any $C^\infty$ family in $\cA_{LF}(\cF,G;\Delta)$
 sufficiently close to of $(\rho_\mu)_{\mu \in \Delta}$
 contains an action parameter-equivalent to $\rho$.

\section{Rigidity and deformation of  flows}
\label{sec:flow}

The real line $\RR$ is the simplest Lie group
 among non-compact and connected ones.
Recall that any locally free $\RR$-action is just a smooth flow
 with no stationary points.
In this section, we discuss
 the rigidity of locally free $\RR$-actions as a model case.

\subsection{Parameter rigidity of locally free $\RR$-actions}
Parameter rigidity of a locally free $\RR$-action
 is characterized by the solvability of a partial differential equation.
\begin{thm}
\label{thm:cohomology flow} 
Let $\rho_0$ be a smooth locally free $\RR$-action
 on a closed manifold $M$
 and $X_0$ the vector field generating $\rho_0$.
Then, $\rho_0$ is parameter rigid
 if and only if the equation
\begin{equation}
\label{eqn:cohomology flow} 
f=X_0g +c.
\end{equation}
 admits a solution $(g,c) \in C^\infty(M,\RR) \times \RR$
 for any given $f \in C^\infty(M,\RR)$.
\end{thm}
The above equation is called {\it the cohomology equation} over $\rho_0$.
\begin{proof}
First, we suppose that $\rho_0$ is parameter-rigid.
Let $\cF$  be the orbit foliation of $\rho_0$
 and take $f \in C^\infty(M,\RR)$.
Since $M$ is closed,
 $f_1=f+c_1$ is a positive valued function
 for some $c_1>0$.
Let $\rho$ be a flow generated by
 a vector field $X=(1/f_1) X_0$.
By the assumption,
 there exists $h \in \Diff_0(\cF)$ and $c_2 \in \RR$ such that
 $\rho^{c_2t} \circ h = h \circ \rho_0^t$.
The diffeomorphism $h$ has the form
 $h(x)=\rho^{-g(x)}$ with some $g\in C^\infty(M,\RR)$.
So, we have
\begin{equation*}
 \rho_0^t(x)=\rho^{c_2t+g \circ \rho_0^t(x)-g(x)}(x)
\end{equation*}
 for any $x \in M$, and hence, $X_0 = (c_2 +X_0 g)X$.
Since $X_0=f_1 X$,
 this implies $f_1=c_2+X_0 g$.
Therefore, the pair $(g,c_2-c_1)$ is a solution of
 (\ref{eqn:cohomology flow}).

Next, we suppose Equation (\ref{eqn:cohomology flow})
 can be solved for any function $f$.
Take an action $\rho \in \cA_{LF}(\cF,\RR)$.
Let $X$ be the vector field generating $\rho$
 and $f$ be the non-zero function satisfying
 $f \cdot X =X_0$.
By assumption, Equation (\ref{eqn:cohomology flow})
 has a solution $(g,c)$ for $f$.
Since $f$ is non-zero and $X_0g(x)=0$ for some $x \in M$, we have $c \neq 0$.
Put $h(x)=\rho^{-g(x)}$.
Then, we have $\rho^{ct} \circ h=h \circ \rho_0^t$.
Since the maps $t \mapsto \rho_0^t(x)$ and
 $t \mapsto \rho^{ct}(h(x))=h(\rho_0^t(x))$ are covering maps
 from $\RR$ to $\cF(x)$ for any $x \in M$, 
 the map $h$ is a self-covering of $M$.
Since $h$ is homotopic to the identity,
 the map $h$ is a diffeomorphism.
Therefore, $\rho$ is parameter equivalent to $\rho_0$.
\end{proof}

We say that a point $x \in M$ is {\it a periodic point}
 of a locally free flow $\rho \in \cA_{LF}(M,\RR)$
 if $\rho^T(x)=x$ for some $T >0$.
The orbit of $x$ is called {\it a periodic orbit}.
A point $x$ is periodic if and only if the orbit $\cO(x,\rho)$ is compact.
\begin{cor}
Let $\rho$ be an action in $\cA_{LF}(M,\RR)$.
Suppose that
 $\rho$ admits two distinct periodic orbits.
Then, $\rho$ is not parameter rigid.
\end{cor}
\begin{proof}
By the assumption, there exists $x_1,x_2 \in M$
 and $T_1,T_2>0$ such that
 $\cO(x_1,\rho) \neq \cO(x_2,\rho)$
 and $\rho^{T_i}(x_i)=x_i$ for each $i=1,2$.
Choose a smooth function $f$
 such that $f \equiv 0$ on $\cO(x_1,\rho)$
 and $f \equiv 1$ on $\cO(x_2,\rho)$.
Then, there exists no solution of 
 (\ref{eqn:cohomology flow}) for $f$.
In fact, if $(g,c)$ is a solution, then we have
\begin{equation*}
 \frac{1}{T}\int_0^T f \circ \rho^t(x) dt=c
\end{equation*}
 for any $x \in M$ and $T >0$ with $\rho^T(x)=x$.
However, the left-hand side should be $0$ or $1$
 for $x=x_1$ or $x_2$.
\end{proof}

There is a classical example of a parameter rigid flow.
For $N \geq 1$, we denote the $N$-dimensional torus
 $\RR^N/\ZZ^N$ by $\TT^N$.
For $v \in \RR^N$,
 we define {\it a linear flow} $R_v$ on $\TT^N$ by $R_v^t(x)=x+tv$.
The vector field $X_v$ corresponding to $R_v$
 is a parallel vector field on $\TT^N$.

We say that $v \in \RR^N$ is {\it Diophantine} if
 there exists $\tau>0$ such that
\begin{equation*}
 \inf_{m \in \ZZ^N \--\{0\}}|\langle m,v \rangle| \cdot \|m\|^\tau>0,
\end{equation*}
 where $\langle,\rangle$
 and $\|\cdot\|$ are the Euclidean inner product 
 and norm on $\RR^N$.
When $v$ is Diophantine,
 we call the flow $R_v$ {\it a Diophantine linear flow}
 and its orbit foliation {\it a Diophantine linear foliation}.

\begin{thm}
[Kolmogorov] 
\label{thm:Kolmogorov}
The cohomology equation (\ref{eqn:cohomology flow})
 over a Diophantine linear flow on $\TT^N$
 admits a solution for any $f \in C^\infty(\TT^N,\RR)$.
By Theorem \ref{thm:cohomology flow},
 any Diophantine linear flow is parameter rigid.
\end{thm} 
\begin{proof}
Take the Fourier expansion
\begin{equation*}
 f(x)=\sum_{m \in \ZZ^N} a_m \exp(2\pi \langle m,x \rangle \sqrt{-1})
\end{equation*} 
 of $f$.
Since $f$ is a smooth function,
 we have
\begin{equation}
\label{eqn:Kolmogorov 1}
 \sup_{m \in \ZZ^N}\|m\|^k |a_m|<\infty
\end{equation}
 for any $k \geq 1$.

Fix a Diophantine vector $v \in \RR^N$.
Put $b_0=0$ and
\begin{equation*}
b_m=\frac{a_m}{2\pi \langle m,v \rangle \sqrt{-1}}
\end{equation*}
 for $m \neq 0$. 
Then,
\begin{equation*}
 g(x)= \sum_{m \in \ZZ^N} b_m \exp(2\pi \langle m,x \rangle \sqrt{-1})
\end{equation*}
 is a formal solution of $f=X_v g + a_0$.
Since $v$ is Diophantine,
 there exists $\tau>0$ and $C>0$ such that
 $|b_m| \leq C\|m\|^\tau |a_m|$
 for any $m \in \ZZ^N$.
By Equation (\ref{eqn:Kolmogorov 1}), we have
\begin{equation*}
 \sup_{m \in \ZZ^N}\|m\|^k |b_m|<\infty
\end{equation*}
 for any $k \geq 1$.
It implies that $g$ is a smooth function.
\end{proof}

Diophantine linear flows are the only known examples
 of parameter rigid flows.
\begin{conj}
[Katok]
Any parameter rigid flow on a closed manifold
 is conjugate to a Diophantine linear flow.
\end{conj}

Recently, some partial results on the conjecture are obtained.
\begin{thm}
[F.Rodrigues-Hertz and A.Rodrigues-Hertz \cite{RHRH06}] 
Let $M$ be a closed manifold with the first Betti number $b_1$.
If $\rho \in \cA_{LF}(M,\RR)$ is parameter rigid,
 then there exists a smooth submersion $\pi:M \ra \TT^{b_1}$
 and a Diophantine linear flow $R_v$ on $\TT^{b_1}$
 such that $\pi \circ \rho^t=R_v^t \circ \pi$.

In particular,
 if $b_1=\dim M$, then $M$ is diffeomorphic to $\TT^{b_1}$
 and $\rho$ is conjugate to a Diophantine linear flow.
\end{thm}

\begin{thm}
[Forni \cite{Fo08}, Kocsard \cite{Ko09}, and Matsumoto \cite{Ma10}] 
Any locally free parameter rigid flow on 
 a three-dimensional closed manifold
 is conjugate to a Diophantine linear flow on $\TT^3$.
\end{thm}

\subsection{Deformation of flows}
\label{sec:flow deform}
There is no known example of a locally rigid flow
 and it is almost hopeless to find it.
\begin{prop}
If $\rho \in \cA_{LF}(M,\RR)$ is locally rigid,
 then there exists a neighborhood $\cU$ of $\rho$
 such that any $\rho' \in \cU$ admits no periodic point.
\end{prop}
\begin{proof}
Let $\cU$ be the conjugacy class of $\rho$.
Since $\rho$ is locally rigid, it is a neighborhood of $\rho$.
For $\rho' \in \cA_{LF}(M,\RR)$, put
\begin{equation*}
\Lambda(\rho')=\{\det D\rho^T_x \st x \in M, T \in \RR, \rho^T(x)=x\}.
\end{equation*}
It is invariant under conjugacy.
Hence, $\Lambda(\rho')=\Lambda(\rho)$ for any $\rho' \in \cU$.

By the Kupka-Smale theorem (see {\it e.g.}, \cite{Ro99}),
 the set $\cU$ contains a flow with at most countably many periodic orbits.
The local rigidity of $\rho$
 implies that $\rho$ admits at most countably many periodic orbits.
Hence, $\Lambda(\rho)$ is at most countable.
However, if $\Lambda(\rho)$ is non-empty,
 then small perturbation on a small neighborhood of a periodic orbit
 can produce a flow $\rho' \in \cU$
 such that $\Lambda(\rho') \neq \Lambda(\rho)$.
\end{proof}

It is unknown whether
 any open subset of $\cA_{LF}(M,\RR)$ contains
 a flow with a periodic point or not.
On the other hand,
 any open subset of the set of $C^1$ flows
 ({\it with $C^1$-topology}) contains a $C^\infty$ flow with a periodic point.
It is just an immediate consequence of Pugh's $C^1$ closing lemma
 \cite{Pu67}.
The validity of the $C^\infty$ closing lemma is
 a long-standing open problem in the theory of dynamical systems.

The following exercise shows that
 it is hard to find a locally complete deformation of a flow.
\begin{excs}
Suppose that a flow $\rho \in \cA_{LF}(M,\RR)$
 admits infinitely many periodic orbits.
Show that any deformation $(\rho_{\mu})_{\mu \in \Delta}$ of $\rho$
 is not locally complete.
\end{excs}

On the other hand, 
 the Diophantine linear flow admits
 a locally transverse deformation.
\begin{thm}
\label{thm:transverse flow}
Let $v \in \RR^N$ be a Diophantine vector
 and $E \subset \RR^N$ be its orthogonal complement.
Then, the deformation $(R_{v+\mu})_{\mu \in E}$ of $R_v$
 is locally transverse.
\end{thm}
Remark that the above deformation is not complete.
In fact, it is easy to see that
 $R_v$ can be approximated by a flow with
 finitely many periodic orbits,
 and hence, which is not conjugate to any linear flow.

The theorem is derived from the following result due to Herman.
Fix $N \geq 2$ and a point $x_0 \in \TT^N$.
Let $\Diff(\TT^N,x_0)$ be the set of diffeomorphisms of $\TT^N$
 which fix $x_0$.
\begin{thm}
[Herman] 
\label{thm:Herman}
Suppose  that $v \in \RR^N$ is Diophantine.
Then, there exists a neighborhood $\cU$ of $X_v$
 in $\frX(\TT^N)$,
 a neighborhood $\cal V$ of $\Id_{\TT^N}$ in $\Diff(\TT^N,x_0)$,
 and a continuous map $\bar{w}:\cU \ra \RR^N$ 
 which satisfy the following property:
 For any $Y \in \cU$,
 there exists a unique diffeomorphism $h \in \cal V$
 such that $Y=h_*(X_v)+X_{\bar{w}(Y)}$.
\end{thm}
\begin{proof}
We give only a sketch of proof here.
See {\it e.g.,} \cite{AG07} for details.
We define a map
 $\Phi:\Diff(\TT^N,x_0) \times \RR^N \ra \frX(\TT^N)$
 by $(h,w) \mapsto h_*(X_v)+X_w$.
The theorem is an immediate consequence of
 the Nash-Moser inverse function theorem
 if we can apply it to $\Phi$ at $(h,w)=(\Id_{\TT^N},v)$.
By the solvability of Equation (\ref{eqn:cohomology flow}) for any $f$,
 we can show that ``the differential'' $D\Phi$ is invertible
 on a neighborhood of $(\Id_{\TT^N},v)$\footnote{Since $\Diff(\TT^N,x_0)$
 is a Fr\'echet manifold (not a Banach manifold),
 the definition of the differential $D\Phi$ is non-trivial.}
 and the inverse satisfies a ``tame'' estimate.
It allows us to apply the Nash-Moser inverse function theorem.
\end{proof}

\begin{proof}
[Proof of Theorem \ref{thm:transverse flow}] 
Let $\cU$, ${\cal V}$, and $\bar{w}$ be the neighborhoods
 and the map in Herman's theorem.
For $\rho \in \cA(\TT^N,\RR)$,
 we denote the vector field generating $\rho$ by $Y_\rho$.
Take neighborhoods $U$ of $0$ in $E$
 and $\cal W$ of a deformation $(R_{v+\mu})_{\mu \in E}$
 in $\cA_{LF}(\TT^N,\RR;E)$, and a constant $\delta>0$
 such that $(1+c)Y_{\rho_\mu} \in \cU$
 for any $(\rho_\mu)_{\mu \in E} \in {\cal W}$,
 $\mu \in U$, and $c \in (-\delta,\delta)$.
For $(\rho_\mu)_{\mu \in E} \in {\cal W}$,
 we define a map $\Psi_{(\rho_\mu)}: U \times (-\delta,\delta) \ra \RR^N$
 by $\Psi_{(\rho_\mu)}(\mu,c)=\bar{w}((1+c) \cdot Y_{\rho_\mu})$.
It is a continuous map which depends continuously on $(\rho_\mu)_{\mu \in E}$.
By the uniqueness of the choice of $h \in \Diff_0(\TT^N)$
 in Herman's theorem,
 we have $\Psi_{(R_{v+\mu})}(\mu,c)=(1+c)\mu+cv$.
In particular, $\Psi_{(R_{v+\mu})}$ is a local homeomorphism
 between neighborhoods of
 $(\mu,c)=(0,0) \in E \times \RR$ and $0 \in \RR^n$.
By the continuous dependence of $\Psi_{(\rho_\mu)}$ with respect to
 the family $(\rho_\mu)$,
 if $(\rho_\mu)_{\mu \in E}$ is sufficiently close to $(R_{v+\mu})_{\mu \in E}$
 then the image of $\Psi_{(\rho_\mu)}$ contains $0$.
In other words, there exists $(\mu_*,c_*) \in U \times (-\delta,\delta)$
 such that $\Psi_{(\rho_\mu)}(\mu_*,c_*)=0$.
Hence, there exists $h_* \in {\cal V}$
 which conjugates $R_v$ with $\rho_{\mu_*}$.
\end{proof}

The above family $(R_{v+\mu})_{\mu \in E}$
 is the best possible in the following sense.
\begin{excs}
Let $(\rho_\mu)_{\mu \in \Delta} \in \cA_{LF}(\TT^N,\RR;\Delta)$
 be a deformation of $R_v$ for $v \in \RR^N$.
Show that if the dimension of $\Delta$ is less than $N-1$,
 then $(\rho_\mu)_{\mu \in \Delta}$ is not a locally transverse deformation.
\end{excs}

\section{The leafwise cohomology}
\label{sec:leafwise cohomology}

As we saw in the previous section,
 the cohomology equation
 plays an important role in the rigidity problem
 of locally free $\RR$-actions.
For actions of abelian Lie groups,
 the solvability of the equation
 is generalized to the almost vanishing of
 the first leafwise cohomology of the orbit foliation.
In this section, we give the definition of the leafwise cohomology
 and show some of its basic properties.
We also compute the cohomology for several examples.

\subsection{The definition and some basic properties}
\label{sec:cohomology}
Let $\cF$ be a foliation on a manifold $M$.
As before, we denote the tangent bundle of $\cF$ by $T\cF$.
We also denote the dual bundle of $T\cF$ by $T^*\cF$.
For $k \geq 0$,
 let $\Omega^k(\cF)$ be the space of smooth sections
 of $\wedge^k T^*\cF$.
Each element of $\Omega^*(\cF)$
 is called {\it a leafwise $k$-form}.

By Frobenius' theorem,
 if $X,Y \in \frX(\cF)$, then $[X,Y] \in \frX(\cF)$.
Hence, we can define {\it the leafwise exterior derivative}
 $d_\cF^k: \Omega^k(\cF) \ra \Omega^{k+1}(\cF)$ by
\begin{align*}
 (d_\cF^k \omega)(X_0,\cdots,X_k)
 & =\sum_{0 \leq i \leq k}
  (-1)^i X_i\left(\omega(X_0,\cdots,\check{X_i},\cdots,X_k)\right)\\
 & +\sum_{0 \leq i<j \leq k}
 (-1)^{i+j}
 \omega([X_i,X_j],X_0,\cdots,\check{X_i},\cdots,\check{X_j},\cdots,X_k)
\end{align*}
 for $X_0,\cdots,X_k \in \frX(\cF)$.
Same as the usual exterior derivative,
 the leafwise derivative satisfies
 $d_\cF^{k+1} \circ d_\cF^k=0$.
For $k \geq 0$,
 {\it the  $k$-th leafwise cohomology group} $H^k(\cF)$
 is the $k$-th cohomology group of
 the cochain complex $(\Omega^*(\cF), d_\cF)$.

\begin{expl}
\label{expl:H0}
$H^0(\cF)$ is the space of smooth functions
 which are constant on each leaf of $\cF$.
Hence, if $\cF$ has a dense leaf, then
 $H^0(\cF) \simeq \RR$.
\end{expl}

\begin{expl}
Suppose that $\cF$ is a one-dimensional orientable foliation on
 a closed manifold $M$.
Let  $X_0$ be a vector field generating $\cF$.
Take $\omega_0 \in \Omega^1(\cF)$ such that
 $\omega_0(X_0) = 1$.
Then, we have $d_\cF^0g=(X_0g) \cdot \omega_0$
 for any $g \in \Omega^0(\cF)=C^\infty(M,\RR)$.
Since $d_\cF^1$ is the zero map,
 the cohomology equation (\ref{eqn:cohomology flow})
 is solved for any $f \in C^\infty(M,\RR)$
 if and only if $H^1(\cF) \simeq \RR$.
In this case, $[\omega_0]$ is a generator of $H^1(\cF)$.
\end{expl}

There are two important homomorphisms whose target is $H^*(\cF)$.
The first is a homomorphism from the de Rham cohomology group.
Let $\Omega^k(M)$ and $H^k(M)$ be the space of (usual) smooth $k$-forms
 and the $k$-th de Rham cohomology group of $M$.
By Frobenius' theorem,
 the restriction of a closed ({\it resp.} exact) $k$-form
 to $\otimes^k T\cF$ defines a $d_\cF$-closed ({\it} resp. exact)
 leafwise $k$-form.
So, the restriction map $r:\Omega^k(M) \ra \Omega^k(\cF)$
 induces a homomorphism $r_*:H^k(M) \ra H^k(\cF)$.

The second is a homomorphism from the cohomology
 of a Lie algebra when $\cF$ is the orbit foliation
 of a locally free action.
Let us recall the definition of the cohomology group
 of a Lie algebra.
Let $\frg$ be a Lie algebra.
For $k \geq 0$, we define the differential
 $d_\frg^k:\wedge^k \frg^* \ra \wedge^{k+1}\frg^*$ by
 $d_\frg^0=0$ and
\begin{equation*}
 (d_\frg^k \alpha)(\xi_0,\cdots,\xi_k)
 =\sum_{0 \leq i<j \leq k}(-1)^{i+j} \alpha([\xi_i,\xi_j],
 \xi_0, \cdots, \check{\xi}_i,\cdots, \check{\xi}_j, \cdots, \xi_k)
\end{equation*}
 for $k \geq 1$ and $\xi_0,\cdots,\xi_k \in \frg$.
The $k$-th cohomology group $H^k(\frg)$
 is the $k$-th cohomology group of the chain complex
 $(\wedge^* \frg^*,d_\frg)$.
\begin{excs}
 $H^1(\frg)$ is isomorphic to $\frg/[\frg,\frg]$.
\end{excs}

Suppose that $\cF$ is the orbit foliation
 of a locally free action $\rho$ of a Lie group $G$.
Let $\frg$ be the Lie algebra of $G$
 and $I_\rho \in \Hom(\frg, \frX(M))$
 be the infinitesimal action associated with $\rho$.
Then, $I_\rho$ induces a homomorphism
 $\iota_\rho: \wedge^* \frg^* \ra \Omega^*(\cF)$ by
\begin{equation*}
 \iota_\rho(\alpha)_x(I_\rho(\xi_1),\cdots I_\rho(\xi_k))
 =\alpha(\xi_1,\cdots, \xi_k)
\end{equation*}
 for any $\alpha \in \wedge^k\frg^*$, $\xi_1,\cdots,\xi_k \in \frg$,
 and $x \in M$.
Since the map $\iota_\rho$ commutes with the differentials,
 it induces a homomorphism $(\iota_\rho)_*:H^*(\frg) \ra H^*(\cF)$.

\begin{prop}
\label{prop:iota}
 The homomorphism $(\iota_\rho)_*:H^1(\frg) \ra H^1(\cF)$
 between the first cohomology groups
 is injective when $M$ is a closed manifold.
\end{prop}
\begin{proof}
Fix $\alpha \in \Ker d_\frg^1$ such that $(\iota_\rho)_*([\alpha]) = 0$.
Then, there exists $g \in C^\infty(M,\RR)$ such that
 $\iota_\rho(\alpha)=d_\cF g$.
For $\xi \in \frg$, let $\Phi_\xi$ be the flow on $M$
 generated by $I_\rho(\xi)$.
For any $\xi \in \frg$, $T>0$, and $x \in M$,
\begin{equation*}
 \alpha(\xi) \cdot T
 =\int_{\{\Phi_\xi^t(x)\}_{0 \leq t \leq T}} \iota_\rho(\alpha)
 =g \circ \Phi_\xi^T(x)-g(x).
\end{equation*}
Since the last term is bounded and $T$ is arbitrary,
 we have $\alpha(\xi)=0$ for any $\xi \in \frg$.
Therefore, $\alpha=0$.
\end{proof}
The same conclusion holds for the homomorphism
 between higher cohomology groups if the action preserves
 a Borel probability measure. See \cite{dS02}.

\begin{expl}
 Let $\cF$ be the orbit foliation of
 a Diophantine linear flow $R_v$ on $\TT^N$.
By Theorem \ref{thm:Kolmogorov},
 $H^1(\cF)$ is isomorphic to $\RR$.
The above proposition implies 
 that $H^1(\cF)$ is generated by
 the dual $\omega_v$ of the constant vector field $X_v$.
The form $\omega_v$ is the restriction of a usual $1$-form.
So, $H^1(\cF)=\Img \iota_*=\Img r_*$.
In particular, the map $r_*$ is not injective for $N \geq 2$.
\end{expl}

The vanishing of the first leafwise cohomology group of the orbit foliation
 implies the existence of an invariant volume.
\begin{prop}
[dos Santos \cite{dS07}] 
\label{prop:inv vol}
Let $G$ be a simply connected Lie group and
 $\cF$ a foliation on an orientable closed manifold $M$.
If $H^1(\cF) \simeq H^1(\frg)$,
 then any $\rho \in \cA_{LF}(\cF,G)$
 preserves a smooth volume,
 {\it i.e.} there exists a smooth volume $\nu$ on $M$
 such that $(\rho^g)^*\nu=\nu$ for any $g \in G$.
\end{prop}
\begin{proof}
Fix an action $\rho \in \cA_{LF}(\cF,G)$
 and a smooth volume form $\nu$ on $M$.
We define a leafwise one-form $\omega \in \Omega^1(\cF)$ by
 ${\cal L}_{X}\nu =\omega(X) \cdot \nu$ for any $X \in \frX(\cF)$.
Then, 
\begin{align*}
(d_{\cF}\omega(X,Y))\cdot \nu
 &= \{X \cdot \omega(Y)-Y \cdot \omega(X)- \omega([X,Y])\}\nu\\
 & = {\cal L}_X ({\cal L}_Y\nu)
   -{\cal L}_Y({\cal L}_X\nu) - {\cal L}_{[X,Y]} \nu\\
 & = 0
\end{align*}
for any $X,Y \in \frX(\cF)$.
Since $H^1(\cF) = \Img (\iota_\rho)_*$
 by assumption and Proposition \ref{prop:iota},
 there exists a smooth function $f$ on $M$
 and $\alpha \in \frg^*$ such that
 $\omega=\iota_\rho(\alpha)+d_\cF f$.
Define a new volume form $\nu_f$ on $M$ by $\nu_f=e^{-f} \cdot \nu$.
It satisfies
\begin{equation*}
 ({\cal L}_{I_\rho(\xi)})\nu_f=\iota_\rho(\alpha)(I_\rho(\xi)) \cdot \nu_f
 =\alpha(\xi) \cdot \nu_f
\end{equation*}
 for any $\xi \in \frg$.
Since $M$ is a closed manifold, $\alpha(\xi)$ must be zero.
It implies that $\rho$ preserves the volume $\nu_f$.
\end{proof}

Remark that the converse of the proposition does not hold.
In fact, there is an easy counterexample.
The linear flow associated with a rational vector
 preserves the standard volume of the torus.
However, the first leafwise cohomology of the orbit foliation
 is infinite dimensional since all points of the torus are periodic.

\subsection{Computation by Fourier analysis}
\label{sec:diophantine}
Theorem \ref{thm:Kolmogorov} can be generalized
 to linear foliations of tori.
Let $B=(v_1,\cdots,v_p) \in \RR^{pN}$ be
 a $p$-tuple of linearly independent vectors in $\RR^N$.
We define {\it the linear action} $\rho_B \in \cA(\TT^N,\RR^p)$
 by $\rho_B^{(t_1,\cdots,t_p)}(x)=x+\sum_{i=1}^pt_iv_i$.
We say that $B=(v_1,\dots,v_p)$ is {\it Diophantine}
 if there exists $\tau>0$ such that
\begin{equation*}
 \inf_{m \in \ZZ^N \--\{0\}}
 \left(\max\{|\langle m,v_1 \rangle|, \cdots, |\langle m,v_p \rangle|\}
 \cdot \|m\|\right)^\tau>0.
\end{equation*}
If $B=(v_1,\dots,v_p)$ is Diophantine,
 the orbit foliation of $\rho_B$ is called
 {\it a Diophantine linear foliation}.
\begin{thm}
[Arraut and dos Santos \cite{AS91}, see also \cite{AS88,ET86}]
\label{thm:linear fol}
Let  $\cF$ be a $p$-dimensional Diophantine linear foliation on $\TT^N$.
Then, $H^1(\cF) \simeq \RR^p$.
\end{thm}
\begin{proof}
Let $B=(v_1,\dots,v_p)$ be a $p$-tuple of linearly independent vectors
 in $\RR^N$ which is Diophantine and whose orbit foliation is $\cF$.
For each $m \in \ZZ^N \--\{0\}$, take $i(m) \in \{1,\dots,p\}$
 such that
\begin{equation*}
|\langle m,v_{i(m)} \rangle| = \max
\{|\langle m,v_1 \rangle|, \dots, |\langle m,v_p \rangle|\}.
\end{equation*}
Since $(v_1,\dots,v_p)$ is Diophantine, there exists $\tau>0$ such that
\begin{equation}
\label{eqn:linear fol 1}
 \inf_{m \in \ZZ^N \--\{0\}} |\langle m,v_{i(m)} \rangle| \cdot \|m\|^\tau>0.
\end{equation}
In particular, $\langle m,v_{i(m)} \rangle \neq 0$
 for any $m \in \ZZ^N \-- \{0\}$.
 
Let $Y_1,\dots,Y_p$ be linear vector fields corresponding to
 $v_1,\dots,v_p$, respectively, and
 $dy_1,\dots,dy_p$ be the dual $1$-forms in $\Omega^1(\cF)$.
Take a closed leafwise $1$-form
 $\omega=\sum_{i=1}^p f_i dy_i$ in $\Omega^1(\cF)$.
Let
\begin{equation*}
f_i(x)=\sum_{m \in \ZZ^N}a_{i,m}\exp(2\pi\langle m,x \rangle \sqrt{-1}) 
\end{equation*}
 be the Fourier expansion of $f_i$.
Since $f_i$ is smooth, we have
\begin{equation}
\label{eqn:linear fol 2}
 \sup_{m \in \ZZ^N}|a_{i(m),m}| \cdot \|m\|^k<+\infty
\end{equation}
 for any $k \geq 1$.
Since $\omega$ is $d_\cF$-closed, we also have $Y_if_j=Y_j f_i$,
 and hence,
\begin{equation}
\label{eqn:linear fol 3}
\langle m,v_{i(m)} \rangle \cdot a_{i,m}
 = \langle m,v_i \rangle \cdot a_{i(m),m}
\end{equation}
 for any $i=1,\dots,p$ and $m \in \ZZ^N \--\{0\}$.
Put $b_0=0$ and $b_m=a_{i(m),m}/\langle m, v_{i(m)} \rangle$
 for $m \in \ZZ^N \--\{0\}$.
By the inequalities (\ref{eqn:linear fol 1}) and (\ref{eqn:linear fol 2}),
\begin{equation*}
 \sup_{m \in \ZZ^N}|b_m| \cdot \|m\|^k<+\infty
\end{equation*}
 for any $k \geq 1$.
Hence, a function
\begin{equation*}
\beta(x)=\sum_{m \in \ZZ^N}b_m\exp(2\pi\langle m,x \rangle \sqrt{-1}) 
\end{equation*}
 is well-defined and smooth.
Since $a_{i,m}=b_m \langle m,v_i \rangle$
 by Equation (\ref{eqn:linear fol 3}),
 we have
\begin{align*}
d_\cF \beta
 & = \sum_{i=1}^p\left(\sum_{m \in \ZZ^N \--\{0\}}b_m \langle m,v_i \rangle
 \exp(2\pi\langle m,x \rangle \sqrt{-1}) \right) dy_i\\
 & = \omega -\sum_{i=1}^p a_{i,0} dy_i.
\end{align*}
Hence, $H^1(\cF)=\Img \iota_{\rho_B}^* \simeq \RR^p$.
\end{proof}

Arraut and dos Santos also computed the higher leafwise cohomology 
 of Diophantine linear foliations.
\begin{thm}
[Arraut and dos Santos \cite{AS91}]
\label{thm:linear fol 2}
Let  $\cF$ be a $p$-dimensional Diophantine linear foliation on $\TT^N$.
Then, $H^*(\cF) \simeq H^*(\TT^p)$.
\end{thm}

The Fourier expansion can be regarded as the irreducible
 decomposition of the regular representation of $\TT^N$.
There is another example of an application of representation theory
 to computation of the leafwise cohomology of a foliation.
Let $\Gamma$ be a cocompact lattice of $SL(2,\RR)$ and put
\begin{equation*}
u(t)=\begin{pmatrix} 1 & t \\ 0 & 1 \end{pmatrix} 
\end{equation*}
 for $t \in \RR$.
We define an action
 $\rho \in \cA_{LF}(\Gamma \bs SL(2,\RR),\RR)$ 
 by $\rho(\Gamma x,t)=\Gamma(x u(t))$.
The $\RR$-action $\rho$ is called {\it the horocycle flow}.
Flaminio and Forni \cite{FF03} gave a detailed description of
 the solution of cohomology equation over the horocycle flow
 by using irreducible decomposition of the regular right
 representation of $SL(2,\RR)$ on $\Gamma \bs SL(2,\RR)$.
When we replace $\RR$ by $\CC$, we obtain a $\CC$-action
 on $\Gamma \bs SL(2,\CC)$.
Its orbit foliation is called {\it the horospherical foliation}.
Using the result by Flaminio and Forni,
 Mieczkowski computed the first cohomology of
 the horospherical foliation.
\begin{thm}
[Mieczkowski \cite{Mi07}]
Let $\cF$ be the orbit foliation of the above $\CC$-action
 on $\Gamma \bs SL(2,\CC)$.
Then, the image of $d_\cF^0$ is a closed subspace of $\Ker d_\cF^1$
 and there exists a subspace $H$ of $\Ker d_\cF^1$ such that
 $H \simeq H^1(M)$ and
\begin{equation*}
 \Ker d_\cF^1 = \Img d_\cF^0 \oplus \Img \iota_{\rho} \oplus H.
\end{equation*}
In particular, $H^1(\cF) \simeq \RR^2 \oplus H^1(M)$.\footnote{
Mieczkowski stated only the isomorphism
$H^1(\cF) \simeq \RR^2 \oplus H^1(M)$
 in \cite[Theorem 1]{Mi07}.
However, by a careful reading of his proof,
 we can see that $\rho$-invariant distributions
 on $C^\infty(\Gamma \bs SL(2,\CC))$
 give the projections associated with
 the splitting $\Ker d_\cF^1=\Img d_\cF^0 \oplus \Img \iota_{\rho} \oplus H$.
}
\end{thm}

\subsection{Computation by the Mayer-Vietoris argument}
\label{sec:calculation}
Let $\cF$ be a foliation on a manifold $M$.
By $\cF|_U$, we denote the restriction of $\cF$ to
 an open subset $U$ of $M$.
More precisely,
 the leaf $(\cF|_U)(x)$ is the connected component
 of $\cF(x) \cap U$ which contains $x$.
For $k \geq 0$, we define a sheaf $\Omega_\cF^k$
 and a pre-sheaf $H_\cF^k$
 by $\Omega_\cF^k(U)=\Omega^k(\cF|_U)$ and $H_\cF^k(U)=H^k(\cF|_U)$.
For open subsets $U_1$ and $U_2$ of $M$,
 we can show the Mayer-Vietoris exact sequence
\begin{equation*}
 \cdots \stackrel{\delta^*}{\ra}
 H^k_{\cF}(U_1 \cup U_2) \stackrel{j^*}{\ra}
  H^k_{\cF}(U_1) \oplus H^k_{\cF}(U_2) \stackrel{i^*}{\ra}
 H^k_{\cF}(U_1 \cap U_2) \stackrel{\delta^*}{\ra}
 H^{k+1}_{\cF}(U_1 \cup U_2) \stackrel{j^*}{\ra} \cdots
\end{equation*}
 same as the de Rham cohomology.

Let us compute the leafwise cohomology of a foliation of suspension type,
 using the above exact sequence.
Let  $\cF$ be a foliation on a manifold $M$.
Suppose that a diffeomorphism $h$ of $M$ satisfies
 $h(\cF(x))=\cF(h(x))$ for any $x \in M$.
By $M_h$, we denote the mapping torus
 $M \times \RR /(x,t+1) \sim (h(x),t)$.
The product foliation $\cF \times \RR$ on $M \times \RR$
 induces a foliation $\cF_h$ on $M_h$.
The foliation $\cF_h$ is called {\it the suspension foliation} of $\cF$.

Take an open cover $M_h=U_1 \cup U_2$
 such that $U_1=M \times (0,1)$ and $U_2= M \times (-1/2,1/2)$.
Then, the natural projection from $U_i$ to $M$
 induces an isomorphism between $H^*(\cF)$ and $H^*_{U_i}(\cF_h)$.
Similarly, $H^*_{U_1 \cap U_2}(\cF_h)$ is
 naturally isomorphic to $H^*(\cF) \oplus H^*(\cF)$.
Under these identifications,
 the map $i^*$ is given by $i^*(a,b)=(a-b,a-h^*(b))$
 for $(a,b) \in H^*(\cF) \oplus H^*(\cF)$.
Hence, we have
\begin{equation*}
\Ker i^* \simeq \Ker (I-h^*),\hsp 
\Img i^* \simeq H^*(\cF) \oplus \Img (I-h^*).
\end{equation*}
The Mayer-Vietoris exact sequence implies
\begin{align*}
H^*(\cF_h)
 & \simeq \Ker i^* \oplus \Img \delta^{*-1}
 \simeq \Ker i^* \oplus [H^{*-1}(\cF)\oplus H^{*-1}(\cF)]/\Img i^{*-1}.
\end{align*}
Therefore,
\begin{equation}
\label{eqn:MV}
H^*(\cF_h) \simeq  \Ker (I-h^*) \oplus [H^{(*-1)}(\cF)/\Img (I-h^{(*-1)})].
\end{equation}

We compute $H^1(\cF_h)$ for the suspension
 of the stable foliation of a hyperbolic toral automorphism.
Let $A$ be an integer valued matrix with $\det A=1$.
We define a diffeomorphism $F_A$ on $\TT^2$ by
 $F_A(x+\ZZ^2)=Ax+\ZZ^2$.
Suppose that eigenvalues $\lambda,\lambda^{-1}$ of $A$
 satisfy $\lambda>1>\lambda^{-1}>0$.
Let $E^s$ be the eigenspace of $\lambda^{-1}$
 and $\cF^s$ be the foliation on $\TT^2$
 given by $\cF^s(x)=x+E^s$.
Since $F_A(\cF^s(x))=\cF^s(F_A(x))$,
 the foliation $\cF^s$ induces
 the suspension foliation $\cF_A$
 on the mapping torus $M_A$.
\begin{thm}
[El Kacimi-Alaoui and Tihami \cite{ET86}] 
\label{thm:ET86}
$H^1(\cF_A) \simeq \RR$.
\end{thm}
\begin{proof}
It is known that $\cF^s$ is a Diophantine linear foliation.
So, we have $H^0(\cF^s) \simeq H^1(\cF^s) \simeq \RR$.
By a direct computation,
 we can check that $F_A^*=I$ on $H^0(\cF^s)$
 and $F_A^*=\lambda^{-1} \cdot I$ on $H^1(\cF^s)$.
The isomorphism (\ref{eqn:MV}) implies $H^1(\cF^s) \simeq \RR$.
\end{proof}
El Kacimi-Alaoui and Tihami also computed
 the first leafwise cohomology group
 for the suspension foliation of higher dimensional
 hyperbolic toral automorphisms.
See \cite{ET86}.

Same as the de Rham cohomology,
 the Mayer-Vietoris sequence for the leafwise cohomology is generalized
 to a spectral sequence.
Let $\cU=\{U_i\}$ be a locally finite open cover of $M$.
By the same construction as the \v{C}ech-de Rham complex
 (see {\it e.g.,} \cite{BT}),
 we obtain a double complex $(C^*(\cU,\Omega_\cF^*),d_\cF,\delta)$,
 where
\begin{align*}
C^p(\cU,\Omega_\cF^q) &
 =\oplus_{i_1<\cdots<i_p}\Omega_\cF^q(U_{i_1} \cap \cdots U_{i_p})
\end{align*}
 and $\delta :C^*(\cU,\Omega_\cF^*) \ra  C^{*+1}(\cU,\Omega_\cF^*)$
 is a natural linear map induced by inclusions.
Moreover, we can show that the sequence
\begin{equation*}
0 \ra \Omega^q(\cF) \ra C^0(\cU, \Omega_\cF^q) 
 \stackrel{\del}{\ra} C^1(\cU,\Omega_\cF^q)
 \stackrel{\del}{\ra} C^2(\cU,\Omega_\cF^q) \ra \cdots
\end{equation*}
 is exact.
The following theorem is proved by the standard method.
\begin{thm}
[El Kacimi-Alaoui and Tihami \cite{ET86}] 
There exists a spectral sequence $\{E_r^{*,*}\}$ such that 
 $E_1^{p,q}=C^p(\cU,H_\cF^q)$,
 $E_2^{p,q}=H_\delta^p(\cU,H_\cF^q)$, and
 $\{E_r^{*,*}\}$ converges to $H^*(\cF)$.
\end{thm}
The reader can find several applications of the spectral sequence
 in \cite{ET86}.

\subsection{Other examples}
\label{sec:calculation2}
In this subsection,
 we give several examples of foliations
 whose first leafwise cohomology is computed by other methods.
\medskip

Fix $p \geq 1$
 and a cocompact lattice $\Gamma$ of $SL(p+1,\RR)$.
Put $M_\Gamma=\Gamma \bs SL(p+1,\RR)$.
By $A$, we denote the subset of $SL(p+1,\RR)$
 that consists of positive diagonal matrices.
It is a closed subgroup of $SL(p+1,\RR)$
 isomorphic to $\RR^p$.
{\it The Weyl chamber flow}
 is the action $\rho \in \cA_{LF}(M_\Gamma, A)$
 given by $\rho(\Gamma x,a)=\Gamma(x a)$.
Let $\cA_p$ be the orbit foliation of $\rho$.
\begin{thm}
[Katok and Spatzier \cite{KS94}] 
\label{thm:KS}
If $p \geq 2$,
 $H^1(\cA_p) \simeq \RR^p$.
\end{thm}
The key features of the proof
 are the decay of matrix coefficient of
 the regular representation 
 and the hyperbolicity of $A$-action.
Remark that Katok and Spatizer proved a similar result
 for a wide class of Lie groups of real-rank more than one.

As an application of the above theorem,
 we compute the first cohomology of another foliation on $M_\Gamma$.
Let $P$ be the subgroup of $SL(p+1,\RR)$ that
 consisting of upper triangular matrices
 with positive diagonals.
It naturally acts on $M_\Gamma$ from right.
Let $\cF_p$ be the orbit foliation of this action.
\begin{thm}
If $p \geq 2$,  $H^1(\cF_p) \simeq \RR^p$.
\end{thm}
\begin{proof}
Let $E_{ij}$ be the square matrix of size $(p+1)$
 whose $(i,j)$-entry is one and the other entries are zero.
For $i,j=1,\cdots,p+1$,
 we define flows $\Phi_{ij}$ and $\Psi_{ij}$ on $M_\Gamma$ by
 $\Phi_{ij}^t(\Gamma g)=\Gamma g \exp(t(E_{ii}-E_{jj}))$
 and $\Psi_{ij}^t(\Gamma g) =\Gamma g \exp(t E_{ij}))$.
Let $X_{ij}$ and $Y_{ij}$ be the vector fields on $M_\Gamma$
 which correspond to $\Phi_{ij}$ and $\Psi_{ij}$, respectively.
Remark that $\cA_p$ is generated by $X_{ij}$'s
 and $\cF_p$ is generated by $X_{ij}$'s and $Y_{ij}$'s with $i<j$.

Take a $d_{\cF_p}$-closed $1$-form $\omega \in \Omega^1(\cF_p)$.
The restriction of $\omega$ to $T\cA_p$ is $d_{\cA_p}$-closed.
By Theorem \ref{thm:KS},
 there exists $h \in C^\infty(M_\Gamma,\RR)$ such that
 $(\omega+dh)(X_{ij})$ is a constant function
 for any $i,j=1,\cdots,p+1$.

We put $\omega'=\omega+dh$ and show $\omega'(Y_{ij})=0$
 for any $i<j$.
Fix $i,j,k=1,\cdots,p+1$ so that $i<j$ and $k \neq i,j$.
Since $[X_{ik},Y_{ij}]=Y_{ij}$, $d_{\cF_p}\omega'(X_{ik},Y_{ij})=0$,
 and $\omega'(X_{ik})$ is constant,
 we have $X_{ik}(\omega'(Y_{ij}))=\omega'(Y_{ij})$.
It implies that 
 $\omega'(Y_{ij})(\Phi_{ik}^t(x)) = e^t \cdot \omega'(Y_{ij})(x)$
 for any $t \in \RR$ and $x \in M_\Gamma$.
By the compactness of $M_\Gamma$, 
 $\omega'(Y_{ij})$ is constantly zero.
Therefore, any $d_{\cF_p}$-closed $1$-form is
 cohomologous to the constant form which vanishes at $Y_{ij}$
 for any $i<j$.
\end{proof}

For $p=1$, the Weyl chamber flow is an $\RR$-action,
 and it is naturally identified with
 the geodesic flow of a two-dimensional hyperbolic orbifold.
It admits infinitely many periodic points,
 and hence, $H^1(\cA_1)$ is infinite dimensional.
By contrast,
 the following theorem asserts that $H^1(\cF_1)$ is finite dimensional.
Let $\rho_\Gamma$ be the natural right action of $P$
 on $\Gamma \bs SL(2,\RR)$.
We denote the Lie algebra of $SL(2,\RR)$ by $\mathfrak{sl}_2(\RR)$.
Let $\iota_{\rho_\Gamma}:H^1(\mathfrak{sl}_2(\RR)) \ra H^1(\cF_\Gamma)$
 and $r_*:H^1(M) \ra H^1(\cF_\Gamma)$ be the homomorphisms
 defined in Section \ref{sec:cohomology}.
\begin{thm}
[Matsumoto and Mitsumatsu \cite{MM03}] 
\label{thm:MM03}
The map
\begin{equation*}
(\iota_{\rho_\Gamma})_* \oplus r_*:
 H^1(\mathfrak{sl}_2(\RR)) \oplus H^1(\Gamma \bs SL(2,\RR))
 \ra H^1(\cF_1) 
\end{equation*}
 is an isomorphism.
\end{thm}
Kanai \cite{Ka09} proved the corresponding result for general
 simple Lie groups of real-rank one.
In the both results,
 the key feature of the proof is
 the hyperbolicity of the $A$-subaction.

\section{Parameter deformation}
\label{sec:orbit preserving}

Now, we come back to the study of
 the deformation of locally free actions.
In this section,
 we discuss parameter rigidity and
 existence of locally complete orbit-preserving deformations.
In the case of abelian actions,
 parameter rigidity is completely characterized by
 the (almost) vanishing of the first cohomology of the orbit foliation.
It was firstly shown by Arraut and dos Santos \cite{AS88}
 for linear foliations on tori,
 and proved by Matsumoto and Mitsumatsu \cite{MM03} for general
 abelian actions.
In Sections \ref{sec:canonical} and \ref{sec:parameter deform},
 we see this characterization.
Because of non-linearity,
 the relationship between parameter rigidity
 and the vanishing of the leafwise cohomology
 is not clear for non-abelian action in general.
However, for several solvable actions,
 vanishing of the cohomology implies parameter rigidity.
We investigate such examples in Sections \ref{sec:non-abelian}
 and \ref{sec:SL}.

\subsection{The canonical one-form}
\label{sec:canonical}
Let $G$ be a simply connected Lie group
 and $\frg$ be its Lie algebra.
To simplify the presentation,
 we assume that $G$ is linear, {\it i.e.},
 a closed subgroup of $GL(N,\RR)$ with some large $N \geq 1$.
Then, each element of $\frg$ is naturally identified
 with a square matrix of size $N$.

Fix a foliation $\cF$ on a closed manifold $M$.
A $\frg$-valued leafwise $1$-form
 $\omega \in \Omega^1(\cF) \otimes \frg$
 is called {\it non-singular} if $\omega_x:T_x\cF \ra \frg$
 is a linear isomorphism for any $x \in M$.
For any action $\rho \in \cA_{LF}(\cF,G)$,
 the associated infinitesimal action $I_\rho$ is non-singular, {\it i.e.},
 the map $(I_\rho)_x:\frg \ra T_x\cF$ is an isomorphism.
Hence, it induces a non-singular $1$-form $\omega_\rho \in \Omega^1(\cF)
 \otimes \frg$ by $(\omega_\rho)_x=(I_\rho)_x^{-1}$.
We call $\omega_\rho$ {\it the canonical $1$-form} of $\rho$.

\begin{lemma}
\label{lemma:canonical form}
Let $\xi_1,\cdots,\xi_p$ be a basis of $\frg$
 and $\alpha_1, \cdots, \alpha_p$ be its dual basis of $\frg^*$.
Then, we have
\begin{equation*}
\omega_\rho=\sum_{i=1}^p \iota_\rho(\alpha_i) \otimes \xi_i, 
\end{equation*}
 where $\iota_\rho:\frg^* \ra \Omega^1(\cF)$
 is the homomorphism defined in Section \ref{sec:cohomology}.
\end{lemma}
\begin{proof}
For $\xi=\sum_{i=1}^p c_i \xi_i$,
 we have $\omega_\rho((I_\rho)_x(\xi))=\xi$
 by definition.
On the other hand,
\begin{equation*}
 \sum_{i=1}^p \iota_\rho(\alpha_i)((I_\rho)_x(\xi)) \otimes \xi_i
 =\sum_{i=1}^p c_i\iota_\rho(\alpha_i)((I_\rho)_x(\xi_i)) \otimes \xi_i
 =\sum_{i=1}^p c_i \xi_i =\xi.
\end{equation*}
\end{proof}

The group of automorphisms of $G$ acts
 (from left) on $\cA_{LF}(\cF,G)$
 by $(\Theta \cdot \rho)(x,g)=\rho(x,\Theta^{-1}(g))$.
\begin{excs}
Show that $\omega_{\Theta \cdot \rho}=\Theta_* \omega_\rho$,
 where $\Theta_*:\frg \ra \frg$ is the differential of $\Theta$.
\end{excs}

The following proposition characterizes the canonical $1$-form.
\begin{prop}
\label{prop:canonical 1}
A $\frg$-valued leafwise $1$-form $\omega \in \Omega^1(\cF) \otimes \frg$
 is the canonical $1$-form of some action in $\cA_{LF}(\cF,G)$
 if and only if it is a non-singular $1$-form which
 satisfies the equation
\begin{equation}
\label{eqn:canonical eq 1} 
d_\cF \omega +[\omega,\omega]=0,
\end{equation}
 where $[\omega,\omega]$ is a $\frg$-valued leafwise two-form
 defined by $[\omega,\omega]_x(v,w)=[\omega(v),\omega(w)]$ 
 for $v,w \in T_x\cF$.
\end{prop}
\begin{proof}
Fix a basis $\xi_1,\cdots,\xi_l$ of $\frg$.
Let $\{c_{ij}^k\}$ be the structure constants of $\frg$,
 {\it i.e.}, $[\xi_i,\xi_j]=\sum_{k=1}^l c_{ij}^k \xi_k$.

Take a non-singular $1$-form $\omega \in \Omega^1(\cF) \otimes \frg$.
Let $X_i$ be a nowhere-vanishing vector field in $\frX(\cF)$ given by
 $X_i(x)=\omega_x^{-1}(\xi_i)$.
Then,
\begin{align*}
(d_\cF\omega+[\omega,\omega])(X_i,X_j) &
 = X_i(\omega(X_j))-X_j(\omega(X_i))
 -\omega([X_i,X_j])+[\omega(X_i),\omega(X_j)]\\
 & = -\omega([X_i,X_j])+\sum_k c_{ij}^k \xi_k\\
 & = -\omega([X_i,X_j])+\sum_k c_{ij}^k \omega(X_k)\\
 &= \omega(\sum_k c_{ij}^k X_k-[X_i,X_j]).
\end{align*}
Since $\omega$ is non-singular,
 $d_\cF\omega+[\omega,\omega]=0$
 if and only if $[X_i,X_j]=\sum_{k}c_{ij}^k X_k$ for any $i,j$.
The latter condition is equivalent to
 that the linear map $\xi_i \mapsto \sum X_i$
 is a homomorphism between Lie algebras.
Hence, $d_\cF\omega+[\omega,\omega]=0$ if and only if
 there exists $\rho \in \cA_{LF}(\cF,G)$
 such that $I_\rho(\xi_i)=X_i$, equivalently,
 $\omega_\rho(X_i(x))=\xi_i=\omega(X_i(x))$ for any $i$.
\end{proof}

The following proposition
 describes how the canonical $1$-form is transformed
 under parameter-equivalence of actions.
We denote the constant map from $M$ to $\{1_G\}$ by $b_{1_G}$.
\begin{prop}
\label{prop:canonical 2}
An action $\rho \in \cA_{LF}(\cF,G)$ is equivalent to $\rho_0$
 if and only if
 there exists a smooth map $b:M \ra G$
 homotopic to $b_{1_G}$
 and an endomorphism $\Theta:G \ra G$ such that
\begin{equation}
\label{eqn:canonical eq 0}
\omega_\rho=b^{-1} \cdot \Theta_*\omega_{\rho_0} \cdot b +b^{-1}d_\cF b. 
\end{equation}
\end{prop} 

To prove the proposition, we need to introduce cocycles over an action.
Let $H$ be another Lie group and $\frh$ be its Lie algebra.
Fix an action $\rho_0 \in \cA_{LF}(\cF,G)$.
We say that $a \in C^\infty(M \times G, H)$ is
 {\it a ($H$-valued) cocycle} over $\rho_0$
 if $a(x,1_G)=1_H$
 and $a(x,gg')=a(x,g) \cdot a(\rho_0^g(x),g')$
 for any $x \in M$ and $g,g' \in G$.
For a cocycle $a$,
 we define {\it the canonical $1$-form}
 $\omega_a \in \Omega^1(\cF) \otimes \frh$
 of $a$ by
\begin{equation*}
(\omega_a)_x(X)=
 \frac{d}{dt}a(x,\exp t(\omega_{\rho_0})_x(X))|_{t=0}.
\end{equation*}

\begin{lemma}
\label{lemma:cocycle coincide}
Two cocycles $a_1$ and $a_2$ over $\rho_0$ coincide
 if $\omega_{a_1}=\omega_{a_2}$.
\end{lemma}
\begin{proof}
For $i=1,2$, we define $\Phi_i:M \times H \times G \ra M \times H$
 by $\Phi_i((x,h),g)=(\rho_0^g(x),h \cdot a(x,g))$.
It is easy to see that $\Phi_i$ is a locally free action and
\begin{equation*}
I_{\Phi_i}(\xi)(x,h)=(I_{\rho_0}(x), h \cdot \omega_{a_i}(I_{\rho_0}(\xi)(x)))
 \in T_xM \times h\cdot \frh \simeq T_{(x,g)}(M \times H).
\end{equation*}
If $\omega_{a_1}=\omega_{a_2}$, then $I_{\Phi_1}=I_{\Phi_2}$,
 and hence, $\Phi_1=\Phi_2$.
It implies $a_1=a_2$.
\end{proof}

Each action in $\cA_{LF}(\cF,G)$ defines
 a $G$-valued cocycle naturally
 and its canonical $1$-form is the canonical $1$-form of the action.
\begin{lemma}
[Arraut and dos Santos \cite{AS90}]
\label{lemma:G-cocycle}
For any $\rho \in \cA_{LF}(\cF,G)$,
 there exists a unique $G$-valued cocycle $a$ over $\rho_0$
 which satisfies $\rho(x,a(x,g))=\rho_0(x,g)$
 for any $x \in M$ and $g \in G$.
Moreover, $a(x,\cdot):G \ra G$ is a diffeomorphism 
 for any $x \in M$.
\end{lemma}
\begin{proof}
For any $x \in M$,
 the maps $\rho_0(x,\cdot), \rho(x,\cdot):G \ra \cF(x)$
 are coverings with $\rho_0(x,1_G)=\rho(x,1_G)=x$.
Since $G$ is simply connected,
 there exists a unique diffeomorphism $a_x$ of $G$
 such that $\rho(x,a_x(g))=\rho_0(x,g)$ and $a_x(1_G)=1_G$.
Put $a(x,g)=a_x(g)$.
It is easy to see that the map $a$ is smooth and satisfies
 $\rho(x,a(x,gg'))=\rho(x,a(x,g) \cdot a(\rho_0(x,g),g'))$.
By the uniqueness of $a_x$, we have
 $a(x,gg')=a(x,g) \cdot a(\rho_0(x,g),g')$,
 and hence, $a$ is a cocycle.

\end{proof}

\begin{lemma}
Let $\rho$ be an action in $\cA_{LF}(\cF,G)$
 and $a:M \times G \ra G$ the cocycle in Lemma \ref{lemma:G-cocycle}.
Then,  $\omega_a$ is equal to the canonical $1$-form of $\rho$.
\end{lemma}
\begin{proof}
Take the differential of the equation
 $\rho(x,a(x,\exp(t \xi)))=\rho_0(x,\exp(t \xi))$ at $t=0$
 for $\xi \in \frg$.
Then, we have $I_\rho(\omega_a(I_{\rho_0}(\xi)))=I_{\rho_0}(\xi)$.
Since $I_{\rho_0}$ is non-singular and $(\omega_\rho)_x=(I_\rho)_x^{-1}$,
 we have $\omega_a=\omega_\rho$.
\end{proof}

\begin{lemma}
\label{lemma:diffeo} 
Let $\rho_1$ and $\rho_2$ be actions in $\cA_{LF}(\cF,G)$,
 $\Theta$ an endomorphism of $G$,
 and $h$ a $C^\infty$ map which is homotopic to the identity.
If $h(\cF(x)) \subset \cF(x)$ for any $x \in M$
 and $\rho_2^{\Theta(g)} \circ h=h \circ \rho_1^{g}$ for any $g \in G$,
 then $h$ is a diffeomorphism and $\Theta$ is an automorphism.
\end{lemma}
\begin{proof}
Since $h$ is homotopic to the identity, it is surjective.
It implies that $h(\cF(x))=\cF(x)$ for any $x \in M$.
If the differential $\Theta_*:\frg \ra \frg$ is not an automorphism,
 then $h(\cF(x))=\{\rho_2^{\Theta(g)}(h(x)) \st g \in G\}$
 is a strict subset of $\cF(h(x))=\cF(x)$ by Sard's theorem.
Hence, $\Theta_*$ must be an automorphism of $\frg$.
Since $G$ is simply connected, $\Theta$ is an automorphism of $G$.

The maps $\rho_1(x,\cdot)$
 and $h \circ \rho_1(x,\Theta(\cdot)) =\rho_2(h(x),\cdot)$
 are covering maps from $G$ to $\cF(x)$.
It implies that $h$ is a self-covering of $M$.
Since $h$ is homotopic to the identity, $h$ is a diffeomorphism.
\end{proof}

Now, we are ready to prove Proposition \ref{prop:canonical 2}.
\begin{proof}
[Proof of Proposition \ref{prop:canonical 2}]
For an endomorphism $\Theta:G \ra G$ and $b \in C^\infty(M,G)$,
 we define a cocycle $a_{b,\Theta}$ over $\rho_0$
 by $a_{b,\Theta}(x,g)=b(x)^{-1} \cdot \Theta(g) \cdot b(\rho_0^g(x))$.
Let  $\omega_{b,\Theta}$ be its canonical $1$-form.
By a direct calculation, we have
\begin{equation*}
\omega_{b,\Theta}=b^{-1}\Theta_*\omega_{\rho_0}b +b^{-1}d_\cF b. 
\end{equation*}

Suppose that $\rho$ is equivalent to $\rho_0$.
Let $h$ a diffeomorphism in $\Diff_0(\cF)$
 and $\Theta$ an automorphism of $G$
 such that $\rho^{\Theta(g)} \circ h=h \circ \rho_0^g(x)$ for any $g \in G$.
Since $h$ is homotopic to $\Id_M$
 through diffeomorphisms preserving each leaf of $\cF$,
 we can take a smooth map $b:M \ra G$ homotopic to $b_{1_G}$
 such that $h(x)=\rho^{b(x)^{-1}}$.
Then,
 $\rho^{b(x)^{-1}\Theta(g)}(x)=\rho^{b(\rho_0^g(x))^{-1}} \circ \rho_0^g(x)$,
 and hence, $\rho_0(x,g)=\rho(x,a_{b,\Theta}(x,g))$.
By Lemma \ref{lemma:G-cocycle},
 the cocycle $a_\rho$ corresponding to $\rho$
 is equal to $a_{b,\Theta}$.
It implies $\omega_{a_\rho}=\omega_{b,\Theta}$,
 and hence,
$\omega_\rho=b^{-1} \Theta_* \omega_{\rho_0}b +b^{-1} d_\cF b$.

Suppose that the equation (\ref{eqn:canonical eq 0}) holds
 for some $\Theta$ and $b$.
Since $\omega_{\alpha_\rho}=\omega_\rho=\omega_{b,\Theta}$,
 the cocycle $a_\rho$ corresponding to $\rho$ coincides with $a_{b,\Theta}$.
It implies that $\rho(x,a_{b,\Theta}(x,g))=\rho_0(x,g)$.
Put $h(x)=\rho^{b(x)^{-1}}$.
Then, we have $\rho^{\Theta(g)} \circ h=h\circ \rho_0^g$.
By Lemma \ref{lemma:diffeo},
 $\Theta$ is an automorphism and $h$ is a diffeomorphism.
\end{proof}

The above interpretation in terms of leafwise $1$-forms
 can be done for general cocycles.
\begin{prop}
[Matsumoto and Mitsumatu \cite{MM03}]
Let $G,H$ be simply connected Lie groups
 and $\frg,\frh$ be their Lie algebras.
Let $\rho_0$ be a locally free action of $G$
 on a closed manifold $M$
 and $\cF$ be the orbit foliation of $\rho_0$.
\begin{enumerate}
\item A $1$-form $\omega \in \Omega^1(\cF) \otimes \frh$
 is the canonical $1$-form of some $H$-valued cocycle over $\rho_0$
 if and only if $d_\cF \omega +[\omega,\omega]=0$.
\item
Let $a_1,a_2$ be $H$-valued cocycles over $\rho_0$,
 $b:M \ra H$ be a smooth map homotopic to $b_{1_H}$,
 and $\Theta$ be an endomorphism of $H$.
Then, the equation
\begin{equation*}
 a_2(x,g)=b(x)^{-1} \cdot \Theta(a_1(x,g)) \cdot b(\rho_0^g(x))
\end{equation*}
 holds if and only if
\begin{equation*}
\omega_{a_2}=b^{-1} (\Theta_*\omega_{a_1}) b +b^{-1} d_\cF b,
\end{equation*}
 where $\omega_{a_i}$ is the canonical $1$-form of the cocycle $a_i$.
\end{enumerate}
\end{prop}

We can extend Propositions \ref{prop:canonical 1}
 and \ref{prop:canonical 2}
 to the case that $G$ may not be a linear group.
In this case, $b^{-1}(\Theta_*\omega_{\rho_0}) b$ is
 replaced by the adjoint ${\mbox Ad}_{b^{-1}}\Theta_*\omega_{r_0}$,
 and $b^{-1}d_\cF b$ is replaced by the pull-back $b^*\theta_G$
 of {\it the Maurer-Cartan form}
 $\theta_G \in \Omega^1(G) \otimes \frg$,
 where $\theta_G(\xi(x))=\xi$ for any $\xi \in \frg$.

\subsection{Parameter deformation of  $\RR^p$-actions}
\label{sec:parameter deform}
Let $M$ be a closed manifold and $\cF$ a foliation on $M$.
Recall that the first cohomology group of $\RR^p$
 as a Lie algebra is $\RR^p$.
Let $\rho$ be an action in $\cA_{LF}(\cF,\RR^p)$,
 $\iota_\rho:\RR^p \ra \Omega^1(\cF)$ is
 the natural homomorphism induced by $I_\rho$,
 and $\omega_\rho$ be the canonical $1$-form.
Since $\Img (\iota_\rho)_* \simeq \RR^p$,
 Lemma \ref{lemma:canonical form} implies
\begin{equation*}
 \Img (\iota_\rho)_* \otimes \RR^p
 = \{\Theta_* \omega_\rho \st \Theta \text{ is a endomorphism of }G\}.
\end{equation*}

Identify the abelian group $\RR^p$
 and the group of positive diagonal matrices of size $p$
 and apply Propositions \ref{prop:canonical 1}
 and \ref{prop:canonical 2} for $\RR^p$-actions.
Then, we obtain the following correspondence
 between actions in $\cA(\cF,\RR^p)$
 and $\RR^p$-valued leafwise $1$-forms.
\begin{prop}
\label{prop:R-k}
A $\RR^p$-valued leafwise $1$-form
 is the canonical $1$-form of an action in $\cA_{LF}(\cF,\RR^p)$
 if and only if it is non-singular and closed.
Two actions $\rho_1,\rho_2 \in \cA_{LF}(\cF,\RR^p)$
 are parameter-equivalent if and only if 
 the cohomology class $[\omega_{\rho_2}]$
 is contained in $\Img (\iota_{\rho_1})_*$.
\end{prop}
As a corollary,
 we obtain a generalization of Theorem \ref{thm:cohomology flow}.
\begin{thm}
[Matsumoto and Mitsumatsu \cite{MM03}, see also \cite{AS88,KS94,dS93}]
\label{thm:MM}
Let $\rho$ be a locally free $\RR^p$-action on a closed manifold
 and $\cF$ be its orbit foliation.
Then $\rho$ is parameter-rigid
 if and only if $H^1(\cF) \simeq \RR^p$.
\end{thm}

\begin{expl}
 Diophantine linear actions on $\TT^N$
 (Theorem \ref{thm:linear fol}) and
 the Weyl chamber flow (Theorem \ref{thm:KS})
 are parameter rigid.
Mieczkowski's action on $M_\Gamma=\Gamma \bs SL(2,\CC)$ is
 also parameter rigid when $H^1(M_\Gamma)$ is trivial.
\end{expl}

What happens for Mieczkowski's example when $H^1(M_\Gamma)$ is non-trivial?
By the following general criterion,
 there exists a locally complete and locally transverse parameter deformation
 parametrized by an open subset of $H^1(M_\Gamma)$.
\begin{thm}
\label{thm:abel canonical}
Let $\cF$ be a foliation on a closed manifold $M$
 and $\rho$ be an action in $\cA_{LF}(\cF,\RR^p)$.
Suppose that $\Img d_\cF^0$ is closed and
 there exists a finite dimensional subspace $H$ of $\Ker d_\cF^1$
 such that
 $\Ker d_\cF^1 = \Img d_\cF^0 \oplus \Img \iota_\rho \oplus H$.
Then, there exists an open neighborhood $\Delta$
 of $0$ in $H \otimes \RR^p$
 and a locally complete and locally transverse parameter deformation
 $(\rho_\mu)_{\mu \in \Delta} \in \cA(\cF,\RR^p;\Delta)$ of $\rho$.
\end{thm}
\begin{proof}
Let $\omega_\rho$ be the canonical $1$-form of $\rho$
 and $\Delta$ be the set of $1$-forms $\mu \in H \otimes \RR^p$
 such that $\omega_\rho+\mu$ is a non-singular $1$-form.
For each $\mu \in \Delta$, 
 there exists the unique action $\rho_\mu \in \cA_{LF}(\cF,\RR^p)$
 whose canonical $1$-form is $\omega_\rho+\mu$.
The set $\Delta$ is an open neighborhood of $0$
 and the family is a parameter deformation of $\rho$

Let us prove the locally completeness of the deformation.
Let $\pi_H:\Ker d_\cF^1 \ra H$ be the projection
 associated with the splitting
 $\Ker d_\cF^1 = \Img d_\cF^0 \oplus \Img \iota_\rho \oplus H$.
It induces a projection
 $\pi_H^{\otimes p}:\Ker d_\cF^1 \otimes \RR^p \ra H \otimes \RR^p$.
It is continuous and the set
\begin{equation*}
{\cal U}=\{\rho' \in \cA_{LF}(\cF,\RR^p)
 \st \pi_H^{\otimes p}(\omega_{\rho'}-\omega_\rho) \in \Delta\}
\end{equation*}
 is an open subset of $\cA_{LF}(\cF,\RR^p)$.
For $\rho' \in {\cal U}$ with $\pi_H^{\otimes p}(\omega_{\rho'})=\mu$,
 the cohomology class $[\omega_{\rho'}-(\omega_\rho+\mu)]$
 is contained in $\Img (\iota_\rho)_*$.
By Theorem \ref{thm:abel canonical},
 $\rho'$ is parameter-equivalent to $\rho_\mu$.
Therefore, $(\rho_\mu)_{\mu \in \Delta}$ is a locally complete deformation.

Next, we show the local transversality.
If a family $(\rho'_\mu)_{\mu \in \Delta}$ is sufficiently close to
 the original family $(\rho_\mu)_{\mu \in \Delta}$,
 then $\{\pi_H^{\otimes p}(\omega_{\rho'_\mu}) \st \mu \in \Delta\}$
 is a neighborhood of $0$ in $H \otimes \RR^p$.
Hence, $[\omega_{\rho'_{\mu_*}}-\omega_\rho] \in \Img (\iota_\rho)_*$
 for some $\mu_* \in \Delta$.
By Theorem \ref{thm:abel canonical} again,
 $\rho'_{\mu_*}$ is parameter-equivalent to $\rho$.
Therefore,  $(\rho_\mu)_{\mu \in \Delta}$ is a locally transverse deformation.
\end{proof}

\subsection{Parameter rigidity of some non-abelian actions}
\label{sec:non-abelian}
As we saw in the previous subsection,
 the equations
 in Propositions \ref{prop:canonical 1} and \ref{prop:canonical 2}
 are linear equations for $\RR^p$-actions.
For general case,
 the equations are non-linear and
 it is unclear whether an action $\rho$ is parameter rigid or not
 even if we know $H^1(\cF)=\Img (\iota_\rho)_*$ for the orbit foliation $\cF$.
However,
 we can reduce the parameter rigidity to the triviality of $H^1(\cF)$
 for several actions of solvable groups.

The first example is an action of the three-dimensional Heisenberg group
\begin{equation*}
 H=\left\{
\left.
\begin{pmatrix}
 1 & x_1 & x_3 \\ 0 & 1 & x_2 \\ 0 & 0 & 1
\end{pmatrix}
 \;\right|\; x_1,x_2, x_3 \in \RR
\right\}.
\end{equation*}
We denote the Lie algebra of $H$ by $\frh$.
\begin{thm}
[dos Santos \cite{dS07}] 
Let $\rho$ be a locally free $H$-action on a closed manifold $M$.
If the orbit foliation $\cF$ of $\rho$
 satisfies that $H^1(\cF) \simeq H^1(\frh)$,
 then $\rho$ is parameter rigid.
\end{thm}
In \cite{dS07}, dos Santos also proved
 the theorem for higher-dimensional Heisenberg groups
 and constructed examples which satisfy the assumption of the theorem.
\begin{proof}
Let
\begin{equation*}
 \xi_1=\begin{pmatrix}	0 & 1 & 0 \\ 0 & 0 & 0 \\ 0 & 0 & 0  \end{pmatrix},
 \hsp
 \xi_2=\begin{pmatrix}	0 & 0 & 0 \\ 0 & 0 & 1 \\ 0 & 0 & 0  \end{pmatrix},
 \hsp
 \xi_3=\begin{pmatrix}	0 & 0 & 1 \\ 0 & 0 & 0 \\ 0 & 0 & 0  \end{pmatrix}
\end{equation*}
 be a basis of $\frh$
 and $\alpha_1,\alpha_2,\alpha_3$ be its dual basis.
Fix $\rho_0 \in \cA_{LF}(M,H)$
 and put $\eta_i=\iota_{\rho_0}(\alpha_i)$ for each $i$.
Since
\begin{equation*}
 [\xi_1,\xi_2]=\xi_3, \hsp
 [\xi_1,\xi_3]=[\xi_2,\xi_3]=0,
\end{equation*}
 we have the equations
\begin{equation*}
 d_\cF \eta_1=d_\cF \eta_2=0, \hsp
 d_\cF \eta_3 =\eta_2 \wedge \eta_1.
\end{equation*}
In particular, $\Img (\iota_{\rho_0})_* \simeq H^1(\frh)$
 is generated by $[\eta_1]$ and $[\eta_2]$.
For $\omega=\sum_{i=1}^3 \omega_i \otimes \xi_i
 \in \Omega^1(\cF) \otimes \frh$,
 the equation $d_\cF \omega+[\omega,\omega]=0$ is equivalent to
\begin{equation}
\label{eqn:Heisenberg 1}
d_\cF\omega_1 = d_\cF\omega_2 = d_\cF \omega_3+\omega_1 \wedge \omega_2=0.
\end{equation}

Fix $\rho \in \cA_{LF}(\cF,H)$.
Let $\omega_\rho=\sum_{i=1}^3 \omega_i \otimes \xi_i$
 be the canonical $1$-form of $\rho$.
First, we will make $\omega_1$ and $\omega_2$
 into forms of $\Img \iota_{\rho_0}$
 by the transformation of canonical $1$-forms
 described in Proposition \ref{prop:canonical 2}.
Since $d_\cF \omega_1=d_\cF \omega_2=0$ by Equation (\ref{eqn:Heisenberg 1})
 and $H^1(\cF)=\Img (\iota_{\rho_0})_*$ by assumption,
 there exists $b_1,b_2 \in C^\infty(M,\RR)$
 and $(c_{ij})_{i,j=1,2} \in \RR^4$ such that
 $\omega_i=c_{i1}\eta_1 + c_{i2}\eta_2 +d_\cF b_i$ for each $i$.
Put
\begin{equation*}
 b(x)=
\begin{pmatrix}
 1 & b_1(x) & 0 \\ 0 & 1 & b_2(x) \\ 0 & 0 & 1 
\end{pmatrix}.
\end{equation*}
By a direct calculation,
 we can show that
 the form $\omega'=\sum_{i,j=1,2} c_{ij}\eta_j \otimes \xi_i
 + \omega'_3 \otimes \xi_3$
 satisfies
\begin{equation}
\label{eqn:Hei 1}
 b^{-1} \omega' b+b^{-1} d_\cF b =\omega_\rho
\end{equation}
 for a suitable choice of $\omega'_3 \in \Omega^1(\cF) \otimes \frh$.

Next, we will make $\omega'_3$ into a $1$-form in $\Img \iota_{\rho_0}$.
Since $\omega'$ satisfies Equation (\ref{eqn:Heisenberg 1}),
\begin{equation*}
 d_\cF\omega'_3=(c_{11}c_{22}-c_{12}c_{21})\eta_2 \wedge \eta_1
 =(c_{11}c_{22}-c_{12}c_{21}) \cdot d_\cF\eta_3.
\end{equation*}
Hence, $\omega'_3-(c_{11}c_{22}-c_{12}c_{21})\eta_3$ is a closed form.
By assumption again,
 there exists $c'_1,c'_2 \in \RR$ and $b'_3 \in C^\infty(M,\RR)$
 such that
\begin{equation*}
 \omega'_3=c_1' \eta_1+c_2' \eta_2+(c_{11}c_{22}-c_{12}c_{21}) \eta_3+d_\cF b'_3. 
\end{equation*}
Put
\begin{equation*}
 b(x)=
\begin{pmatrix}
 1 & 0 & b'_3(x) \\ 0 & 1 & 0 \\ 0 & 0 & 1 
\end{pmatrix}.
\end{equation*}
 and
\begin{equation*}
\omega''=\sum_{i,j=1}^2c_{ij}\eta_j \otimes \xi_i
 +[c_1' \eta_1+c_2' \eta_2+(c_{11}c_{22}-c_{12}c_{21})\eta_3] \otimes \xi_3.
\end{equation*}
Then, we have
\begin{equation}
\label{eqn:Hei 2}
 b^{-1} \omega'' b+b^{-1}d_\cF b =\omega'.
\end{equation}

Finally, we take an endomorphism $\Theta$ of $\mathfrak h$
 such that $\Theta_*(\xi_j)=c_{1j}\xi_1 +c_{2j}\xi_2+c_j'\xi_3$
 for $j=1,2$.
It satisfies $\Theta_*(\xi_3)
 =\Theta_*[\xi_1,\xi_2]=(c_{11}c_{22}-c_{12}c_{21})\xi_3$.
Hence, 
\begin{equation}
\label{eqn:Hei 3}
 \Theta_*\omega_{\rho_0}
 =\Theta_*(\sum_{j=1}^3 \eta_j \otimes \xi_j)
 =\omega''.
\end{equation}
The equations (\ref{eqn:Hei 1}), (\ref{eqn:Hei 2}),
 and (\ref{eqn:Hei 3}) imply 
\begin{equation*}
 \omega_\rho=(bb')^{-1}\Theta_*\omega_{\rho_0}(bb')+ (bb')^{-1}d_\cF(bb').
\end{equation*}
By Proposition \ref{prop:canonical 2},
 the action $\rho$ is equivalent to $\rho_0$.
\end{proof}

Recently, Maruhashi generalized dos Santos' results
 to general simply connected nilpotent Lie groups.
\begin{thm}
[Maruhashi \cite{Mr-pre}] 
\label{thm:Maruhashi}
Let $G$ be a simply connected nilpotent Lie group
 with the Lie algebra $\frg$
 and $\rho$ a locally free $G$-action on a closed manifold $M$.
If the orbit foliation $\cF$ of $\rho$
 satisfies that $H^1(\cF) \simeq H^1(\frg)$, then $\rho$ is parameter rigid.
The converse is true if $\cF$ has a dense leaf.
\end{thm}
He also gave a family of parameter-rigid actions of nilpotent groups
 by generalizing dos Santos' examples for the Heisenberg groups.
Ram\'irez gave another natural action of a nilpotent Lie group
 which satisfies the above condition.
\begin{thm}
[Ram\'irez \cite{Ra09}]
Let $U$ be the nilpotent subgroup of $SL(n,\RR)$ consisting
 of upper triangular matrices whose diagonal entries are one,
 ${\mathfrak u}$ the Lie algebra of $U$,
 and $\Gamma$ a cocompact lattice of $SL(n,\RR)$.
If $n \geq 4$, then the orbit foliation $\cF$ of
 the natural right $U$-action on $\Gamma \bs SL(n,\RR)$
 satisfies $H^1(\cF) \simeq H^1({\mathfrak u})$.
By Theorem \ref{thm:Maruhashi}, the action is parameter rigid.
\end{thm}

The second example we discuss is
 an action of a two-dimensional solvable group
\begin{equation*}
 GA=\left\{\left.
\begin{pmatrix}
 e^t & u \\ 0 & 1 \end{pmatrix}
 \;\right|\; u,t \in \RR \right\}.
\end{equation*}
Let $A$ be an element of $SL(2,\RR)$
 such that the eigenvalues $\lambda,\lambda^{-1}$ are real and
 $\lambda>1$.
Let $F_A$ be  a diffeomorphism of $\TT^2$
 given  by $F_A(z +\ZZ^2)=Az + \ZZ^2$
 and let $M_A$ be the mapping torus
\begin{equation*}
M_A= \TT^2 \times \RR / (x,s+\log \lambda) \sim (F_A(x), s).
\end{equation*}
We define an action $\rho_A \in \cA_{LF}(M_A,GA)$ by
\begin{equation*}
 \rho_A\left([x,s],
  \begin{pmatrix} e^t & u \\ 0 & 1 \end{pmatrix}
 \right)=[x+(e^s u) \cdot v,s+t], 
\end{equation*}
 where $v$ is the eigenvector associated with $\lambda^{-1}$.
Remark that the orbit foliation $\cF$ of $\rho_A$
 is diffeomorphic to
 the second example in Section \ref{sec:calculation}.
\begin{thm}
[Matsumoto-Mitsumatsu \cite{MM03}]
\label{thm:Sol}
The action $\rho_A$ is parameter rigid.
\end{thm}
\begin{proof}
The Lie algebra $\mathfrak{ga}$ of $GA$ has a basis
\begin{equation*}
 \xi_1=
 \begin{pmatrix}
 1 & 0 \\ 0 & 0
 \end{pmatrix}
\hsp
 \xi_2=
 \begin{pmatrix}
 0 & 1 \\ 0 & 0
 \end{pmatrix}.
\end{equation*}
Let $\alpha_1,\alpha_2$ be the dual basis of $\mathfrak{ga}^*$.
We put $\eta_i=\iota_{\rho_A}(\alpha_i)$.
Then, $[\xi_1,\xi_2]=\xi_2$, and hence,
\begin{equation*}
 d_\cF \eta_1=d_\cF \eta_2 +\eta_1 \wedge \eta_2 =0.
\end{equation*}
In particular, we have $\Img (\iota_{\rho_A})_*=[\eta_1]$.
%

Take $\rho \in \cA_{LF}(\cF,GA)$.
Let $\omega_{\rho_A}$ and $\omega_\rho$
 be the canonical $1$-forms of $\rho_A$ and $\rho$.
Then, $\omega_{\rho_A} =\eta_1 \otimes \xi_1+\eta_2 \otimes \xi_2$
 and $\omega_{\rho} = \omega_1 \otimes \xi_1+ \omega_2 \otimes \xi_2$
 for some $\omega_1,\omega_2 \in \Omega^1(\cF)$.
Since $\omega_\rho$ satisfies
 the equation $d_\cF \omega_\rho +[\omega_\rho,\omega_\rho]=0$,
 the form $\omega_1$ is closed.
By Theorem \ref{thm:ET86}, $H^1(\cF)=\Img (\iota_{\rho_A})_*=\RR [\eta_1]$.
Hence, there exists $c_1 \in \RR$ and $b_1 \in C^\infty(M_A,GA)$
 such that $\omega_1=c_1\eta_1+d_\cF b$.
By Proposition \ref{prop:inv vol},
 $\rho$ preserves a smooth volume naturally.
As a (not immediate) consequence of this fact,
 we can obtain $c_1=1$ (see \cite[p.1863--1864]{MM03} for details).
Put $\omega' = \eta_1 \otimes \xi_1 +e^{b_1}\omega_2 \otimes \xi_2$ and
\begin{align*}
 b & =
\begin{pmatrix}
 e^{b_1} & 0 \\ 0 & 1
\end{pmatrix}.
\end{align*}
Then, by a direct calculation, we have
 $b^{-1}\omega'b +b^{-1}d_\cF b=\omega_1 \otimes \xi_1
 +\omega_2 \otimes \xi_2 =\omega_{\rho}$.
Take $f,g \in C^\infty(M,\RR)$
 such that $e^{b_1}\omega_2=f\eta_1 +g \eta_2$.
Since $d_\cF \omega'+[\omega',\omega']=0$,
 the pair $(f,g)$ satisfies
\begin{equation}
\label{eqn:Sol 2}
 X g=S f, 
\end{equation}
 where $X=I_{\rho_A}(\xi_1)$ and $S=I_{\rho_A}(\xi_2)$.

Let $\Theta$ be an endomorphism of $GA$.
Then, $\Theta_*(\xi_1)=\xi_1$ and $\Theta_*(\xi_2)=c_2 \cdot \xi_2$
 for some $c_2 \in \RR$.
For $b' \in C^\infty(M,GA)$ of the form
\begin{equation*}
 b'(x)=
\begin{pmatrix}
 1 & h(x)-c_2 \\ 0 & 1 
\end{pmatrix},
\end{equation*}
 we have
\begin{equation*}
 (b')^{-1} \Theta_*(\omega_{\rho_A})b'
 + (b')^{-1}d_\cF b'
 = \eta_1 \otimes \xi_1 + [(h+X h)\eta_1+(S h-c_2) \eta_2] \otimes \xi_2.
\end{equation*}
Hence, the equivalence of $\rho$ and $\rho_A$
 is reduced to the solvability of an inhomogeneous linear equation
\begin{equation}
\label{eqn:Sol 3}
\left\{
\begin{array}{l}
f = h + X h \\ g = S h - c_2.
\end{array}
\right.
\end{equation}
In fact,
 the following proposition guarantees the solvability,
 and it completes the proof.
\begin{prop}
[Matsumoto and Mitsumatu \cite{MM03}]
\label{prop:MM03}
If smooth functions $f,g$ satisfy 
 the equation (\ref{eqn:Sol 2}),
 then the equation (\ref{eqn:Sol 3}) has a solution $(h,c_2)$.
\end{prop}
\end{proof}
The group $GA$ is naturally isomorphic to
 the subgroup of $SL(2,\RR)$ which consists of
 upper triangular matrices by the map
\begin{equation}
\label{eqn:Sol SL}
\theta:
\begin{pmatrix}
e^t & u \\ 0 & 1
\end{pmatrix} 
 \mapsto
\begin{pmatrix}
 e^{\frac{t}{2}} & e^{-\frac{t}{2}}u \\ 0 & e^{-\frac{t}{2}}
\end{pmatrix}.
\end{equation}
Let $\Gamma$ be a cocompact lattice of $SL(2,\RR)$
 and put $M_\Gamma=\Gamma \bs SL(2,\RR)$.
We define an action $\rho_\Gamma \in \cA_{LF}(M_\Gamma,GA)$
 by $\rho_\Gamma(\Gamma x, g)=\Gamma (x \cdot \theta(g))$.
It is just the second example in Section \ref{sec:calculation2}.
In \cite{MM03},
 Matsumoto and Mitsumatsu showed
 an analogue of Proposition \ref{prop:MM03} for $\rho_\Gamma$.
\begin{prop}
\label{prop:SL Sol}
Let $\xi_1, \xi_2$ be the basis of $\mathfrak{ga}$
 given in the proof
 of Theorem \ref{thm:Sol}.
Put $X=I_{\rho_\Gamma}(\xi_1)$ and $S=I_{\rho_\Gamma}(\xi_2)$.
Then, if smooth functions $f,g \in C^\infty(M_\Gamma,\RR)$
 satisfy $Sf=Xg$, then the equation
\begin{equation}
\left\{
\begin{array}{l}
f = h + Xh \\ g = Sh + c.
\end{array}
\right.
\end{equation}
 has a solution $(h,c) \in C^\infty(M_\Gamma,\RR) \times \RR$.
\end{prop}
When $H^1(M_\Gamma)$ is trivial,
 we have $H^1(\cF) \simeq \RR$ by Theorem \ref{thm:MM03}.
In this case, we can prove the parameter rigidity
 of $\rho_\Gamma$ by an argument similar to the above.
\begin{thm}
[{\it c.f.}, \cite{Gh85,MM03}]
When $H^1(M_\Gamma)$ is trivial,
 then $\rho_\Gamma$  is parameter rigid.
\end{thm}

\subsection{A complete deformation for actions of $GA$}
\label{sec:SL}

Let $\Gamma$ be a cocompact lattice of $SL(2,\RR)$
 and put $M_\Gamma=\Gamma \bs SL(2,\RR)$.
Let $\rho_\Gamma \in \cA_{LF}(M_\Gamma,GA)$ be the action
 given by $\rho_\Gamma(\Gamma x,g)=\Gamma (x \cdot \theta(g))$,
 which is discussed 
 in the last paragraph of the previous subsection.
It is natural to ask whether $\rho_\Gamma$ is parameter rigid
 or not when $H^1(M_\Gamma)$ is non-trivial.

Let $\cF$ be the orbit foliation of $\rho_\Gamma$.
First, we determine the space of {\it infinitesimal}
 parameter deformations in terms of the leafwise cohomology.
Recall that the space $\cA_{LF}(\cF,GA)$ is identified
 with the solution of the non-linear equation
\begin{equation}
\label{eqn:coho 1}
 d_\cF \omega +[\omega,\omega] =0.
\end{equation}
 in $\Omega^1(\cF) \otimes \mathfrak{ga}$.
Two actions are parameter equivalent with trivial automorphism
 if and only if the equation
\begin{equation}
\label{eqn:coho 0}
 \omega_2 = b^{-1}\omega_ 1b +b^{-1}d_\cF b,
\end{equation}
 admits a smooth solution $b:M_\Gamma \ra GA$,
 where $\omega_1$ and $\omega_2$ are the canonical $1$-forms of actions.
Let $\omega_0$ be the canonical $1$-form of $\rho_\Gamma$.
Put $\omega_t=\omega_0+t\omega$ and $b_t=\exp(t \beta)$
 with $\omega \in \Omega^1(\cF) \otimes \frga$
 and $\beta \in \Omega^0(\cF) \otimes \frga$.
Substitute $\omega_t$ and $b_t$ into the above equations
 and take the first order term with respect to $t$.
Then, we obtain the formally linearized equations
{\addtocounter{equation}{-2}
\renewcommand{\theequation}{\arabic{equation}L}
\begin{gather}
 \omega_2-\omega_1 = [\omega_0,\beta]+d_\cF \beta,\\
 d_\cF\omega +[\omega,\omega_0]+[\omega_0,\omega] =0.
\end{gather}
}
We define the linear map
 $d_{\rho_\Gamma}^k:\Omega^k(\cF) \otimes \frga
 \ra \Omega^{k+1}(\cF) \otimes \frga$
 for $k=1,2$ by
 $d_{\rho_\Gamma}^0\beta=[\omega_0,\beta]+d_\cF \beta$
 and 
 $d_{\rho_\Gamma}^1\omega =d_\cF\omega +[\omega,\omega_0]+[\omega_0,\omega]$.
They satisfy $d_{\rho_\Gamma}^1 \circ d_{\rho_\Gamma}^0=0$
 and the above linearized equations become
 $\omega_2-\omega_2  = d_{\rho_\Gamma}^0 \beta$
 and $d_{\rho_\Gamma}^1 \omega =0$.
We call the quotient space $\Ker d_{\rho_\Gamma}^1 /\Img d_{\rho_\Gamma}^0$
 {\it the space of infinitesimal parameter deformations}
 of $\rho_\Gamma$ and we denote it by $H^1(\rho_\Gamma,\cF)$.

\begin{prop}
$H^1(\rho_\Gamma,\cF) \simeq H^1(M_\Gamma)$.
\end{prop}
\begin{proof}
Fix a basis 
\begin{equation*}
\xi_X=
\begin{pmatrix}
1/2 & 0 \\ 0 & -1/2
\end{pmatrix},
\xi_S=
\begin{pmatrix}
0 & 1 \\ 0 & 0
\end{pmatrix},
\xi_U=
\begin{pmatrix}
0 & 0 \\ 1 & 0
\end{pmatrix}
\end{equation*}
 of the Lie algebra $\mathfrak{sl}(2,\RR)$ of $SL(2,\RR)$.
The standard right $SL(2,\RR)$-action on $M_\Gamma$ induces
 vector fields $X, S$, and $U$ which
 correspond to $\xi_X, \xi_S$, and $\xi_U$.
Let $\eta,\sigma$, $\upsilon$ be the dual $1$-forms of
 $X, S$ and $U$, respectively.
Then, $\Omega^1(\cF)$ is generated by $\eta$ and $\sigma$
 as a $C^\infty(M_\Gamma,\RR)$-module.
If $\omega=\omega_X \otimes \xi_X + \omega_S \otimes \xi_S$ is
 $d_{\rho_\Gamma}$-closed, then $d_\cF \omega_X=0$ and
 $d_\cF\omega_S= -(\omega_x(X)+\omega_S(S))\eta \wedge \sigma$.

First, we claim that $\omega=\omega_X \otimes \xi_X+\omega_S \otimes \xi_S$
 is $d_{\rho_\Gamma}$-exact if and only if $\omega_X$ is $d_\cF$-exact.
For $\varphi \in C^\infty(M_\Gamma,\RR)$, we have
\begin{equation*}
 d_{\rho_\Gamma}^0(\varphi \otimes \xi_X)
 =(d_\cF \varphi) \otimes \xi_X
 + (d_\cF \psi + \psi \cdot \eta -\varphi \cdot \sigma) \otimes \xi_S.
\end{equation*}
Hence, if $\omega$ is $d_{\rho_\Gamma}$-exact then $\omega_X$ is $d_\cF$-exact.
Suppose that $\omega_X$ is $d_\cF$-exact.
Take $\varphi \in C^\infty(M_\Gamma,\RR)$ such that $d_\cF \varphi=\omega_X$.
By replacing $\omega$ with $\omega+d_{\rho_\Gamma}^0 (\varphi \otimes \xi_X)$,
 we may assume that $\omega_X=0$.
Put $\omega_S=f \eta + g \sigma$.
Since $\omega$ is $d_{\rho_\Gamma}$-closed, we have $Sf=Xg$.
Proposition \ref{prop:SL Sol} implies
 that there exists $h \in C^\infty(M_\Gamma,\RR)$ and $c \in \RR$
 such that $f=h+X h$ and $g=S h-c$.
Hence, $\omega=d_{\rho_\Gamma}^0(-c \otimes \xi_X+h \otimes \xi_S)$.
It completes the proof of the claim.

By the claim,
 $H^1(\rho_\Gamma,\cF)$ is isomorphic to
\begin{equation*}
 \{[\omega_X] \in H^1(\cF) \st
 d^1_{\rho_\Gamma}(\omega_X \otimes \xi_X +\omega_S \otimes \xi_S)=0
 \text{ for some } \omega_S \in \Omega^1(\cF)\}.
\end{equation*}
So, it is sufficient to show that
 for any $d_\cF$-closed $1$-form
 $\omega_X \in \Omega^1(\cF)$, there exists $\omega_S \in \Omega^1(\cF)$
 such that $\omega=\omega_X \otimes \xi_X+\omega_S \otimes \xi_S$
 is $d_{\rho_\Gamma}^1$-closed.
Fix a Riemannian metric on $M_\Gamma$ such that
 $(X_\Gamma,(S_\Gamma+U_\Gamma)/2,(S_\Gamma-U_\Gamma/2))$
 is an orthonormal framing of $TM_\Gamma$.
By Theorem \ref{thm:MM03}, there exists $f_0 \in C^\infty(M_\Gamma,\RR)$
 such that $\omega_X+d_\cF f_0$ extends to
 a harmonic $1$-form with respected to the metric.
Replacing $\omega_X$ with $\omega_X+d_{\cF}f_0$,
 we may assume that $\omega_X$ is the restriction
 of a harmonic form $\omega_h$ to $T\cF$.
Put $\omega_h=f \eta +g \sigma+ h\upsilon$.
Since $\omega_h$ is harmonic and $M_\Gamma$ is compact
 (it implies ${\cal L}_{(S_\Gamma-U_\Gamma)}\omega_h=0$),
 we can show the equations $2f =(S-U)g$ and $2Yf = -(S+U)g$.
Now, it is easy to check
 that $d_{\rho_\Gamma}(\omega_X \otimes \xi_X+ (-g\eta+f\sigma) \otimes \xi_S)=0$.
\end{proof}

One may expect the existence of a complete deformation
 whose parameter space of an open subset of
 $H^1(M_\Gamma) \simeq H^1(\rho_\Gamma,\cF)$.
The author of this note proved the existence of {\it globally} complete
 deformation.
\begin{thm}
[Asaoka, in preparation] 
There exists an open subset $\Delta_\Gamma$
 of $H^1(M_\Gamma)$ containing $0$
 and a parameter deformation 
 $(\rho_\mu)_{\mu \in \Delta} \in \cA(M_\Gamma,GA;\Delta_\Gamma)$
 of $\rho_\Gamma$ such that
\begin{enumerate}
\item if $\rho_\mu$ is equivalent to $\rho_\nu$
 then $\mu=\nu$, and
\item any $\rho \in \cA_{LF}(\cF,GA)$  is equivalent to
 $\rho_\mu$ for some $\mu \in \Delta_\Gamma$.
\end{enumerate}
\end{thm}
\begin{cor}
[Asaoka \cite{As09}]
When $H^1(M_\Gamma)$ is non-trivial,
 then $\rho_\Gamma$ is not parameter rigid.
\end{cor}
Construction of the deformation $(\rho_\mu)_{\mu \in \Delta_\Gamma}$
 is essentially done in \cite{As09}.
Remark that the proof does not use
 the computation of $H^1(\rho_\Gamma,\cF)$.
It heavily depends on the ergodic theory of hyperbolic dynamics,
 especially on the existence of the Margulis measure,
 and the deformation theory of low dimensional Anosov systems.
To obtain the smoothness of the family,
 we also use the smooth dependence of the Margulis measure,
  in some sense, with respect to the parameter.

It is natural to expect that
 an analogous result holds for $SL(2,\CC)$.
However, the corresponding action for $SL(2,\CC)$
 is locally parameter rigid.
\begin{thm}
[Asaoka \cite{As-pre}] 
\label{thm:Aspre}
Let $\Gamma$ be a cocompact lattice of $SL(2,\CC)$
 and $GA_\CC$ be the subgroup of $SL(2,\CC)$
 which consists of upper triangular matrices.
Then, the standard $GA_\CC$ action on $\Gamma \bs SL(2,\CC)$
 is locally parameter rigid.
\end{thm}

\section{Deformation of orbits}
\label{sec:general}

In this section,
 we discuss deformations which may not
 preserve the orbit foliation.
The equations we need to solve are non-linear even for $\RR^p$-actions,
 as we investigated deformation of linear flows on tori
 in Section \ref{sec:flow}.
The main techniques to describe such deformation
 are linearization and the Nash-Moser type theorems.
The former reduces the problem
 to computation of the bundle-valued leafwise cohomology.
The latter allows us to construct solutions of
 the original non-linear problem from the linear one.

\subsection{Infinitesimal deformation of foliations}
To know deformations of a given locally free action,
 it is natural to investigate deformations of the orbit foliation.
In this subsection,
 we describe the space of infinitesimal deformations
 of a foliation in terms of the leafwise cohomology.

Let $\cF$ be a foliation on a manifold $M$.
To simplify, we assume that $\cF$ admits
 a {\it complementary foliation} $\cF^\perp$,
 {\it i.e.,} it is transverse to $\cF$
 and satisfies $\dim \cF +\dim \cF^\perp =\dim M$.
The normal bundle $TM/T\cF$ of $T\cF$
 can be naturally identified with
 the tangent bundle $T\cF^\perp$ of $\cF^\perp$.
By $\pi^\perp$, we denote
 the projection from $TM=T\cF \oplus T\cF^\perp$ to $T\cF^\perp$.
Let $\Omega^k(\cF;T\cF^\perp)$ be the space
 of $T\cF^\perp$-valued leafwise $k$-forms.
We define the differential
 $d_\cF^k:\Omega^k(\cF;T\cF^\perp) \ra \Omega^{k+1}(\cF;T\cF^\perp)$
 by
\begin{align*}
(d_\cF^k \omega)(X_0,\cdots,X_k)
 & = \sum_{0 \leq i \leq k}
 (-1)^i \pi^\perp
 \left[X_i,\omega(X_0,\cdots,\check{X}_i,\cdots,X_k)\right]\\
 & \quad +\sum_{0 \leq i<j \leq k}
 (-1)^{i+j} \omega([X_i,X_j],X_0,\cdots,\check{X}_i,\cdots,
 \check{X}_j, \cdots,X_k).
\end{align*}
It satisfies $d_\cF^{k+1} \circ d_\cF^k=0$.
We denote the quotient $\Ker d_\cF^k / \Img d_\cF^{k-1}$
 by $H^k(\cF;T\cF^\perp)$.

Suppose that the foliation $\cF$ is $p$-dimensional.
For $\omega \in \Omega^1(\cF;T\cF^\perp)$,
 we define a $p$-plane field $E_\omega$ on $M$ by
\begin{equation*}
E_\omega(x)=\{v+\omega(v) \st v \in T_x\cF\}.
\end{equation*}
It gives a one-to-one correspondence between
 $T\cF^\perp$-valued leafwise $1$-forms and
 $p$-plane field transverse to $T\cF^\perp$.
By a direct computation on a local coordinate
 adapted to the pair $(\cF,\cF^\perp)$,
 we obtain the following criterion for the integrability of $E_\omega$.
\begin{lemma}
The $p$-plane field $E_\omega$ generates a foliation 
 if and only if $\omega$ satisfies the equation
\begin{equation*}
 d_\cF \omega +[\omega,\omega]=0.
\end{equation*}
\end{lemma}

Fix $\beta \in \frX(\cF^\perp) =\Omega^0(\cF;T\cF^\perp)$.
Let $\{h_t\}_{t \in \RR}$ be a one-parameter family
 of diffeomorphisms such that
 $h_0$ is the identity map
 and $h_t$ preserves each orbit of $\cF^\perp$ for any $t$.
We define a family $\{\omega_t\}_{t \in \RR}$
 of $1$-forms in $\Omega^1(\cF;T\cF^\perp)$
 by $E_{\omega_t}=(h_t)_*(T\cF)$
 and a vector field $\beta \in \Omega^0(\cF;T\cF^\perp)$ 
 by $\beta(x)=(d/dt)h_t(x)|_{t=0}$.
By a direct computation
 on a local coordinate adapted to the pair $(\cF,\cF^\perp)$
 again, we have
\begin{equation*}
\lim_{t \ra 0}\frac{1}{t}\omega_t= d_\cF^0 \beta.
\end{equation*}
So, one can regard the cohomology group $H^1(\cF;T\cF^\perp)$
 as the space of infinitesimal deformation
 of the foliation $\cF$.
We say that a foliation $\cF$ {\it infinitesimally rigid}
 if $H^1(\cF;T\cF^\perp)=\{0\}$.

\begin{expl}
Let $\cF$ be the orbit foliation
 of a Diophantine linear action in $\cA_{LF}(\TT^N,\RR^p)$.
Since $T\cF^\perp$ is a trivial bundle,
 Theorem \ref{thm:linear fol} implies
\begin{equation*}
 H^1(\cF;\cF^\perp) \simeq H^1(\cF) \otimes \RR^{N-p}
 \simeq \RR^{N-p}.
\end{equation*}
In particular, $\cF$ is not infinitesimally rigid.
\end{expl}

\begin{excs}
Let $\cF_A$ be the suspension foliation
 associated to a hyperbolic automorphism on $\TT^2$,
 which is defined in Section \ref{sec:calculation}.
Show that $\cF_A$ is infinitesimally rigid
 using a Mayer-Vietoris argument as in Section \ref{sec:calculation}.
\end{excs}

\begin{expl}
[Kononenko \cite{Ko99}, Kanai \cite{Ka09}]
Let  $\cA_p$ be the orbit foliation of the Weyl chamber flow,
 which is defined in Section \ref{sec:calculation2}.
If $p \geq 2$, then $\cA_p$ is infinitesimally rigid.
\end{expl}

\subsection{Hamilton's criterion for local rigidity}
Let $\cF$ be a foliation on a closed manifold $M$
 and $\cF^\perp$ be its complementary foliation.
We say that $\cF$ is {\it locally rigid}
 if any foliation $\cF'$ sufficiently close to $\cF$
 is diffeomorphic to $\cF$.

Using Hamilton's implicit function theorem for non-linear exact sequence
 \cite[Section 2.6]{Ha77},
 one obtain the following criterion for local rigidity of a foliation.
\begin{thm}
[Hamilton \cite{Ha}] 
Suppose that there exist continuous linear operators
 $\delta^k:\Omega^{k+1}(\cF;T\cF^\perp) \ra \Omega^k(\cF;T\cF^\perp)$
 for $k=1,2$,
 an integer $r \geq 1$,
 and a sequence $\{C_s\}_{s \geq 1}$ of positive real numbers such that
\begin{enumerate}
\item $d_\cF^0 \circ \delta^0+\delta^1 \circ d_\cF^1=\Id$,
\item $\|\delta^0 \omega\|_s \leq C_s\|\omega\|_{s+r}$
 and $\|\delta^1 \sigma\|_s \leq C_s\|\sigma\|_{s+r}$
 for any $s \geq 1$,
 $\omega \in \Omega^1(\cF;T\cF^\perp)$,
 and $\sigma \in \Omega^2(\cF;T\cF^\perp)$,
 where $\|\cdot\|_s$ is the $C^s$-norm
 on $\Omega^k(\cF;T\cF^\perp)$.
\end{enumerate}
Then, $\cF$ is locally rigid.

Moreover, we can choose the diffeomorphism $h$ in the definition
 of locally rigidity so that it is close to the identity map.
\end{thm}
 
\begin{thm}
[El Kacimi Alaoui and Nicolau \cite{EN93}] 
Let $\cF_A$ be the suspension foliation related to
 a hyperbolic toral automorphism, which is given
 in Section \ref{sec:calculation}.
Then,  $\cF_A$ satisfies Hamilton's criterion above.
In particular, it is locally rigid.
\end{thm}
With the parameter rigidity of
 the action $\rho_A$ (Theorem \ref{thm:Sol}), we obtain 
\begin{cor}
[Matsumoto and Mitsumatsu \cite{MM03}]
The action $\rho_A$  is locally rigid.
\end{cor}
In \cite{EN93} and \cite{MM03},
 they also proved the corresponding results
 for higher dimensional hyperbolic toral automorphisms.

It is unknown
 whether the orbit foliation of the Weyl chamber flow
 satisfies Hamilton's criterion or not.
However, 
 Katok and Spatzier proved
 the rigidity of the orbit foliation by another method.
\begin{thm}
[Katok and Spatzier \cite{KS97}] 
The orbit foliation $\cA_p$ of the Weyl chamber flow
 is locally rigid if $p \geq 2$.
\end{thm}
With the parameter rigidity of
 the Weyl chamber flow (Theorem \ref{thm:KS})
 we obtain 
\begin{cor}
The Weyl chamber flow is locally rigid if $p \geq 2$.
\end{cor}

\subsection{Existence of locally transverse deformations}
Although deformation theory is well-developed
 for transversely holomorphic foliations,
 ({\it e.g.}, \cite{DK79,DK80,DK84,EN89,Gi92,GHS83,GN89}),
 there is no general deformation theory for smooth foliations
 with non-trivial infinitesimal deformation so far
 since we can not apply Hamilton's criterion in this case.
However, there are several actions
 for which we can find a locally transverse deformation.
One example is a Diophantine linear flow,
 which we discussed in Section \ref{sec:flow}.
In this subsection, we give two more examples.

The first example is a codimension-one Diophantine linear action.
By $\Diff_0(S^1)$, we denote the set of orientation-preserving
 diffeomorphisms of $S^1$.
Let $\cF$ be a codimension-one foliation on $\TT^{p+1}$
 which is transverse to $\{x\} \times S^1$ for any $x \in \TT^p$.
For each $i=1,\cdots,p$,
 we can define a holonomy map $f_i \in \Diff_0(S^1)$ of $\cF$
 along the $i$-th coordinate.
The family $(f_1, \cdots, f_p)$ is pairwise commuting.
On the other hand, when
 a pairwise commuting family $(f_1,\cdots,f_p)$ in $\Diff_0(S^1)$ 
 is given,
 then the suspension construction gives
 a codimension-one foliation on $\cF$,
 which is transverse to $\{x\} \times S^1$ for any $x \in \TT^p$.
Two foliations are diffeomorphic to each other
 if the corresponding
 families $(f_1,\cdots,f_p)$ and $(g_1,\cdots,g_p)$
 are conjugate,
 {\it i.e.,}
 there exists $h \in \Diff_0(S^1)$ such that
 $g_i \circ h = h \circ f_i$ for any $i=1,\cdots,p$.
So, the local rigidity problem of $\cF$ is reduced
 to the problem on a pairwise commuting family $(f_1,\cdots,f_p)$.

For $f \in \Diff_0(S^1)$,
 {\it the translation number}
 $\tau(f) \in \RR/\ZZ$ is defined by
\begin{equation*}
 \left(\lim_{n \ra \infty} \frac{\tilde{f}^n(0)}{n} \right)+\ZZ,
\end{equation*}
 where $\tilde{f}:\RR \ra \RR$ is a lift of $f$,
It is known that the map $\tau:\Diff_0(S^1) \ra S^1$ is continuous
 (see {\it e.g.,} \cite[Proposition 11.1.6]{KH}).
For $\theta \in S^1$,
 let $r_\theta$ be the rotation defined by $r_\theta(x)=x+\theta$.
\begin{thm}
[Moser \cite{Mo90}] 
Let $(f_1,\cdots,f_p)$ be a pairwise commuting family in $\Diff_0(S^1)$.
Suppose that $(1,\tilde{\tau}(f_1),\dots, \tilde{\tau}(f_p)) \in \RR^{p+1}$
 is a Diophantine vector (see Section \ref{sec:diophantine}
 for the definition),
 where $\tilde{\tau}(f_i) \in \RR$ is a representative of
 $\tau(f_i) \in \RR/\ZZ$.
Then, there exists $h \in \Diff_0(S^1)$ such that
 $f_i \circ h = h \circ r_{\tau(f_i)}$
 for any $i=1,\cdots,p$.
\end{thm}
As a consequence of the theorem,
 we can show the existence of locally transverse deformation
 of a codimension-one Diophantine linear action.
\begin{thm}
Let $\rho$ be the linear action of $\RR^p$
 on $\TT^{p+1}$ determined by linearly independent vectors
 $v_1,\cdots,v_p \in \RR^{p+1}$.
Take $w \in \RR^{p+1}$ so that
 $v_1,\cdots,v_p,w$ is a basis of $\RR^{p+1}$
 and we define a $C^\infty$ family of
 actions $(\rho_s)_{s \in \RR^p} \in \cA_{LF}(\TT^{p+1},\RR^p;\RR^p)$ by
\begin{equation*}
\rho_s^t(x)=x+ \sum_{i=1}^p t_i(v_i+s_iw),
\end{equation*}
 for $x \in M$, $t=(t_1,\cdots,t_p), s=(s_1,\cdots,s_p) \in \RR^p$.
If the linear action $\rho_0$ is Diophantine,
 then $(\rho_s)_{s \in \RR^p}$ is locally transverse at $s=0$.
\end{thm}
\begin{excs}
Prove the theorem.
One way is a modification of
 the proof of Theorem \ref{thm:transverse flow}.
One can prove `local transversality of the orbit foliation'
 by the continuity of the rotation number and Moser's theorem,
 in stead of Herman's theorem.
The local transversality of action
 will follow from the parameter rigidity of a Diophantine linear actions.
\end{excs}

The second example is a $\RR^2$-action on
 $\Gamma \bs SL(2,\RR) \times SL(2,\RR)$
 by commuting parabolic elements.
Put
\begin{equation*}
 u^t = \begin{pmatrix} 1 & t \\ 0 & 1\end{pmatrix},
 \hsp
 u_\mu^t = \exp\left(t
 \begin{pmatrix} 0  & 1 \\ \mu & 0 \end{pmatrix}
 \right).
\end{equation*}
Remark that $u_0^t=u^t$.

Let $\Gamma$ be an irreducible cocompact lattice
 of $SL(2,\RR) \times SL(2,\RR)$
 and put $M_\Gamma=\Gamma \bs (SL(2,\RR) \times SL(2,\RR))$.
For $(\mu,\lambda) \in \RR^2$, we define an action
 $\rho_{\mu,\lambda} \in \cA_{LF}(M_\Gamma,\RR^2)$ by
\begin{equation*}
 \rho_{\mu,\lambda}(\Gamma(x,y),(s,t))
 =\Gamma(x u_\mu^s,\, y u_\lambda^t).
\end{equation*}
Let $\cF$ be the orbit foliation of $\rho_{0,0}$.
\begin{thm}
[Mieczkowski \cite{Mi07}] 
$H^1(\cF) \simeq \RR^2$.
In particular,
 the action $\rho_{0,0}$ is parameter rigid.
\end{thm}
One may wish to prove the local transversality of the deformation
 $(\rho_{\mu,\lambda})_{(\mu,\nu) \in \RR^2}$ of $\rho_{0,0}$
 like Diophantine linear actions.
However, we can not apply the techniques for Diophantine linear actions 
 because of the non-linearity of the space $SL(2,\RR)$.
Damjanovi\'c and Katok developed a new Nash-Moser-type scheme
 and they obtained the local transversality.
\begin{thm}
[Danjanovi\'c and Katok \cite{DKpre}] 
The deformation $(\rho_{\mu,\lambda})_{(\mu,\lambda) \in \RR^2}$
 of $\rho_{0,0}$ is locally transverse.
\end{thm}
In \cite{Mi07} and \cite{DKpre},
 they also show parameter rigidity
 and existence of a transverse deformation for another actions
 using the same method.

\subsection{Transverse geometric structures}
In this subsection,
 we sketch another method to describe deformations
 of an orbit foliation which is not locally rigid.

Fix a torsion-free cocompact lattice $\Gamma$
 of $PSL(2,\RR)=SL(2,\RR)/\{\pm I\}$.
It acts on the hyperbolic plane $\HH^2$ naturally
 and $\Sigma=\Gamma \bs \HH^2$ is a closed surface
 of genus $g \geq 2$.
Let $\cT(\Sigma)$ be the Teichm\"uller space of $\Sigma$
 (see {\it e.g.} \cite[Chapter 4]{Jo}
 for the definition and basic properties).
It can be realized as
 a set of homomorphisms $\mu$ from $\Gamma$ to $PSL(2,\RR)$
 whose image $\Gamma_\mu$ is a cocompact lattice.
It is known that $\cT(\Sigma)$ has a natural structure of
 $(6g-6)$-dimensional smooth manifold.

Let $P$ be the subgroup of $PSL(2,\RR)$
 which consists of upper triangular matrices.
For each $\mu$ in $\cT(\Sigma)$,
 we define an action
 $\rho_\mu \in \cA_{LF}(\Gamma_\mu \bs PSL(2,\RR),P)$ by
 $\rho_\mu(\Gamma_\mu x,p)=\Gamma_\mu(x \cdot p)$.
The action is essentially same one as in Sections
 \ref{sec:non-abelian} and \ref{sec:SL}.
Let $\cF_\mu$ be the orbit foliation of $\rho_\mu$.
To simplify notation,
 we put $\rho_\Gamma=\rho_{\text{Id}_\Gamma}$
 and $\cF_\Gamma=\cF_{\text{Id}_\Gamma}$.

It is well-known that the foliation $\cF_\Gamma$ is not locally rigid.
In fact, $M_{\mu_1}$ is diffeomorphic to $M_{\mu_2}$
 for any $\mu_1,\mu_2 \in \cT(\Sigma)$.
However, $\cF_{\mu_1}$ is diffeomorphic to $\cF_{\mu_2}$
 if and only if $\Gamma_{\mu_1}$ is conjugate to $\Gamma_{\mu_2}$
 as a subgroup of $PSL(2,\RR)$.
Hence, the family $\{\cF_\mu\}_{\mu \in \cT(\Sigma)}$
 gives a non-trivial deformation of $\cF_\Gamma$.
Ghys proved that this is the only possible one.

\begin{thm}
[Ghys \cite{Gh92}] 
\label{thm:Ghys local}
Any two-dimensional foliation on $M_\Gamma$ 
 sufficiently close to $\cF_\Gamma$
 is diffeomorphic to $\cF_\mu$
 for some $\mu \in \cT(\Sigma)$.
\end{thm}
He also proved global rigidity.
\begin{thm}
[Ghys \cite{Gh93}] 
\label{thm:Ghys global}
If a two-dimensional foliation $\cF$ on $M_\Gamma$ 
 has no closed leaves,
 then $\cF$ is diffeomorphic to $\cF_\mu$
 for some $\mu \in \cT(\Sigma)$.
\end{thm}
The orbit foliation of a locally free $P$-action has no closed leaf.
Hence, we obtain
\begin{cor}
For any $\rho \in \cA_{LF}(M_\Gamma,P)$,
 there exists $\mu \in \cT(\Sigma)$ such that
 the orbit foliation of $\rho$
 is diffeomorphic to $\cF_\mu$.
\end{cor}

The basic idea of the proof is
 to find a transverse projective structure of the foliation.
Once it is shown, it is not so hard to show that
 $\cF$ is diffeomorphic to $\cF_\mu$ for some $\mu$.
Ghys constructed the transverse projective structure
 by using the theory of hyperbolic dynamical systems.
Kononenko and Yue \cite{KY97} gave an alternative
 proof of Theorem \ref{thm:Ghys local}\footnote{
They proved the $C^3$ conjugacy of foliations.
However, it must be the $C^\infty$ conjugacy
 by a regularity theorem of conjugacies between Anosov flows
 by de la Llave and Moriy\'on \cite{LM88}.}.
They used the vanishing of a twisted cohomology of the lattice $\Gamma$,
 which is closely related
 to the leafwise cohomology of $\cF_\Gamma$ valued in
 the symmetric two-forms on the normal bundle of $T\cF$.
So, it may possible to reduce Theorem \ref{thm:Ghys local} 
 to the vanishing of the bundle-valued leafwise cohomology.

Modifying the construction of a complete parameter deformation
 of $\rho_\Gamma$ (Theorem \ref{thm:Aspre}),
 we obtain a {\it globally} complete deformation of $\rho_\Gamma$.
\begin{thm}
[Asaoka, in preparation] 
There exists an open subset $\Delta$ of $\cT(\Sigma) \times H^1(M_\Gamma)$
 and a $C^\infty$ family $(\rho_\mu)_{\mu \in \Delta} \in \cA_{LF}(M_\Gamma,P)$
 such that any $\rho \in \cA_{LF}(M_\Gamma,P)$
 is conjugate to $\rho_\mu$ for some $\mu \in \Delta$.
\end{thm}

\end{document}